\documentclass[10,reqno]{amsart}

\usepackage{epsfig}

\usepackage{amssymb}
\let\phi\varphi
\let\eps\varepsilon

\def\om {\omega }
\def\Om {\Omega }
\def\pd {\partial }
\let\si\sigma
\let\Si\Sigma
\def\a {\alpha}
\def\l {\lambda}

\def\fP{\mathfrak{P}}

\def\cD{{\mathcal D}}

\def\cP{{\mathcal P}}
\def\cB{{\mathcal B}}
\def\cL{{\mathcal L}}

\def\bR{{\mathbb R}}
\def\bC{{\mathbb C}}

\def\bP{{\mathbb P}}
\def\bE{{\mathbb E}}

\def\div{\mathop{\rm div}}

\def\Image{{\mathop{\rm Im}}}
\def\Int{\mathop{\rm Int}}
\def\Real{{\mathop{\rm Re}}}

\let\le\leqslant
\let\ge\geqslant
\let\<\langle
\let\>\rangle

\newtheorem{deff}{Definition}[section]

\newtheorem{lem}{Lemma}[section]

\newtheorem{theorem}{Theorem}[section]

 \makeatletter \@addtoreset{equation}{section}
\@addtoreset{theorem}{section}

\begin{document}

 \title[Stochastic feedback stabilization]
 {Feedback stabilization  for Oseen fluid equations: \\A stochastic approach}



\author{Jinqiao Duan}
\address[J. Duan]
{Department of Applied Mathematics\\
 Illinois Institute of Technology\\
   Chicago, IL 60616, USA.
{\em E-mail address,} J. Duan: duan@iit.edu}

\author{Andrei\,V.\,Fursikov }
\address[A.V. Fursikov]
{Department of   Mechanics and Mathematics\\
 Moscow State   University\\
  119899 Moscow, Russia.  }
\email[A.V. Fursikov]{fursikov@mtu-net.ru }

\thanks{This work was partially supported
  by a COBASE grant from the National Research Council and   the NSF grant DMS-0209326.}


\maketitle

\begin{abstract}
{\bf J. Math. Fluid Mechanics, 2004, in press.}

 The authors consider stochastic aspects of the stabilization
problem for two and three-dimensional Oseen equations with help of
feedback control defined on a part of the fluid boundary.
Stochastic issues arise when inevitable unpredictable fluctuations
in numerical realization of stabilization procedures are taken
into account and they are supposed to be independent identically
distributed random variables. Under this assumption the solution
to the stabilization problem obtained via boundary feedback
control can be described by a Markov chain or a discrete random
dynamical system. It is shown that this random dynamical system
possesses a unique, exponentially attracting, invariant measure,
namely, this random dynamical system is ergodic. This gives
adequate statistical description of the stabilization process on
the stage when stabilized solution has to be retained near zero
(i.e. near unstable state of equilibrium).

\end{abstract}

\bigskip


\bigskip










\section{Introduction}

This paper is devoted to study some statistical aspects of the
stabilization problem for two and three-dimensional  Oseen
equations with help of feedback control defined on a part  of the
boundary that restricts a domain where the equations are
determined. Since for stabilization problem the case of unstable
equations is interesting, we assume that Oseen equations possess a
solution that is exponentially growing as time $t\to \infty$, i.e.
the solution is unstable. New approach to the problem of
stabilization by feedback control was proposed by one of the
authors of this paper in \cite{F1}-\cite{F7}.  Namely,
construction of stabilization from a part of boundary for
parabolic equation, 2D Oseen equations as well as for 2D and 3D
Navier-Stokes system was created in \cite{F1}, \cite{F2},
\cite{F4}, \cite{F5}. This construction reduces solution of a
stabilization problem to solving a mixed boundary value problem
defined on an extended domain with zero Dirichlet boundary
conditions and with initial condition belonging to a stable
invariant manifold defined in an neighborhood of steady-state
solution, near which we stabilize our system. In the case of
(linear) Oseen equations, steady state solution equals to zero and
invariant manifold is replaced on subpace $X_\sigma$ invariant
with respect to resolving semigroup and such that the solution
going out initial condition $w_0\in X_\sigma$ tends to zero as
time $t\to \infty$ with the rate $e^{-\sigma t}$ or faster.

Evidently, aforementioned mixed boundary value problem is not
stable because if $w_0\not \in X_\sigma$ then solution
$w(t,\cdot)$ outgoing $w_0$ goes away from $X_\sigma$ (and goes
away from zero) even if $w_0$ is arbitrarily close to $X_\sigma$.
That is why straightforward application of proposed construction
of stabilization to numerical simulation may not be successful,
because  unpredictable fluctuations inevitably arise in real
calculations. This situation was analysed in \cite{F3}, \cite{F6},
and \cite{F7} with help of the concept of real process introduced
there. The method of damping of unpredictable fluctuations by some
feedback mechanism was worked out in these papers and an estimate
of stabilization for real process was obtained.

This estimate is informative only when a norm of stabilized real
process is not too  small. But when this norm has the same order
as norms of unpredictable fluctuations, aforementioned estimate
became uninformative. Actually in this situation behavior of
stabilized real process became chaotic. The goal of this paper is
just to investigate the behavior of stabilized real process in
small neighborhood of zero. More precisely, we solve the problem
of retention of stabilized flow near unstable state of
equilibrium. To do this we impose additional assumption on
unpredictable fluctuations. We suppose that they are independent
identically distributed (i.i.d.) random variables with probability
distribution supported in a small neighborhood of origin for phase
space. Then the real process is described by a random dynamical
system and forms a Markov chain.

Our aim is to prove that this Markov chain is ergodic, i.e. it
possesses unique stationary measure invariant with respect to the
corresponding Markov semigroup.  This gives us the possibility to
calculate by well known formulas the statistical characteristics
of stabilized real process and to make clear its behavior using
these formulas.

During the last few years uniqueness   of invariant measures for
2D Navier-Stokes equations have been proved by F.Flandoli,
B.Maslovski \cite{FM}, S.Kuksin, A.Shirikyan \cite{KSh1},
\cite{KSh}, \cite{KSh2}, W.E, J.Mattingly, Ya.Sinai \cite{EMS},
S.Kuksin \cite{K}, Duan and Goldys \cite{Duan-Goldys}, and   other
authors.

To prove ergodicity of indicated Markov chain we use recent
results of S.Kuksin, A.Shirikyan \cite{KSh} and S.Kuksin \cite{K}
on uniqueness of invariant measures for 2D Navier-Stokes equations
with random kick-forces where coupling approach was applied.

Actually, it is enough for us to verify that random dynamical
system arising in stabilization construction indicated above
satisfies all conditions imposed in \cite{KSh}, \cite{K} on random
dynamical systems.

In section 2 we recall   necessary information on stabilization
method. In section 3 we formulate the main results and present
conditions imposed in \cite{KSh}, \cite{K} on random dynamical
systems in a form convenient for our situation. In sections 4-6 we
verify that these conditions fulfil for RDS corresponding to
stabilization proceedure.

We thank S.B.Kuksin and A.Shirikyan for useful discussion on
ergodicity of Navier-Stokes equations. We thank also
M.S.Agranovich for allowing us to read the proof of his result
announced in \cite{A1} and formulated in Lemma \ref{lemAgr} below
before its publication.


\section{Preliminaries to the stabilization theory}

In this section we recall briefly results of \cite{F1}-\cite{F7}
used below.


\subsection{Formulation of the problem.}

Let $\Omega \subset \mathbb{R}^d,\, d=2,3, \; \partial \Omega \in
C^\infty$, $Q=\mathbb{R}_+\times \Omega$. We consider the Oseen
equations:
\begin{equation}\label{2.1}
  \partial _tv(t,x)-\Delta v+(a(x),\nabla )v+(v,\nabla )a+\nabla p(t,x)=0,\quad
\div v(t,x)=0
\end{equation}
with initial condition
\begin{equation}\label{2.3}
  v(t,x)|_{t=0}=v_0(x).
\end{equation}
Here $(t, x)=(t, x_1,\dots , x_d)\in Q,\; v(t,x)=(v_1,\dots ,
v_d)$ is a velocity of fluid flow, $  p(t,x)$ is   pressure,
$a(x)=(a_1(x),\dots , a_d(x))$ is a given solenoidal vector field.

We suppose that $ \partial \Omega=\overline \Gamma \cup \overline
\Gamma _0 , \; \Gamma \ne \emptyset $ where $\Gamma , \Gamma _0$
are open sets (in topology of $\partial \Omega$). Here, as usual,
the over line   means the closure of a set. We define $\Si
=\bR_+\times \Gamma , \Si _0 =\bR_+\times \Gamma _0$, and   set:
\begin{equation}\label{2.5}
  v|_{\Si _0}=0, \quad v|_{\Si }=u
\end{equation}
where $u$ is a control, supported on~$\Si$.

Let  $\sigma >0$ be given. The problem of stabilization with
 rate $\sigma$ of a solution to problem (\ref{2.1})--(\ref{2.5}) is to
construct a control $u$ defined on $\Si$ such that the solution
$v(t,x)$ of boundary value problem (\ref{2.1})--(\ref{2.5})
satisfies:
\begin{equation}\label{2.6}
  \| v(t,x)\|^2_{L_2(\Omega)}\le ce^{-\sigma t}
\end{equation}
where $c>0$ depends on $v_0, \sigma $ and $\Gamma _0$.


\subsection{The main idea of the stabilization method.}

Let $\omega \subset \bR^d$ be a bounded domain such that $\Omega
\cap \omega =\emptyset, \quad
  \overline \Omega \cap \overline \omega=\overline \Gamma .$
We set
\begin{equation}\label{2.8}
  G= \Int (\overline \Omega \cup \overline \omega)
\end{equation}
(the notation  $\Int A$ means, as always, the interior of the set
$A$).

We suppose that $\partial G\in C^\alpha$ where $\alpha >2$ is
fixed and in all points except $\overline \Gamma \setminus \Gamma
\equiv \partial \Gamma $ it possesses the $C^\infty$ smoothness.
For the construction of $\omega $ and detailed description of
$\partial G$ in a neighborhood of $\partial \Gamma $, please see \cite{F2},
  \cite{F6}.

We extend problem (\ref{2.1})--(\ref{2.3}) from $\Omega$  to $G$.
Let us assume that
\begin{equation}\label{2.9}
  a(x)\in V^2(G)\cap \left( H^1_0(G)\right)^d, \quad a(x) \; \mbox{is real valued},
\end{equation}
where, as usually, $H^k(G ),\; k\in \mathbb{N}$, is the Sobolev
space of scalar functions, defined and square integrable on $G$
together with all its derivatives up to order $k$ and $(H^k (G )
)^d $ is the analogous  space of vector fields. Besides,
$H^1_0(G)=\{ f(x)\in H^1(G): f(x)|_{x\in \partial G}=0\}$ and
$$
V^k(\Omega) = \left\{ v(x) = (v_1,\dots , v_d) \in (H^k (\Omega
))^d :
     \quad {\rm div}\, v = 0 \right\}
$$

The extension of (\ref{2.1})--(\ref{2.3}) from $\Omega$  to $G$
can be written as follows:
\begin{equation}\label{2.10}
  \partial _tw(t,x)-\Delta   w+(a(x),\nabla )w+
(w,\nabla )a+\nabla p(t,x)=0,\quad  \div w(t,x)=0,
\end{equation}
\begin{equation}\label{2.12}
  w(t,x)|_{t=0}=w_0(x),\qquad   w|_S=0,
\end{equation}
where $S=\mathbb{R}_+\times \partial G$. Note that, actually,
$w_0$ from (\ref{2.12}) will be some special extension of $v_0(x)$
from (\ref{2.3}) such that the solution $w(t,x)$ of problem
(\ref{2.10}), (\ref{2.12})  satisfies the inequality
\begin{equation}\label{2.39}
  \|w(t,\cdot) \|_{V^0 (G)}\le ce^{-\sigma t} \|
     w_0 \|_{V^0 (G)}\quad \mbox{for}\ t\ge 0
\end{equation}

For vector fields defined on $G$ we denote by $\gamma _{\Omega }$
the operator of restriction on $\Omega$ and by $\gamma _{\Gamma }$
we denote the operator of restriction on~$\Gamma$:
\begin{equation}\label{2.91}
  \gamma _{\Omega}: V^k(G)\longrightarrow V^k(\Omega), \quad
  \gamma _{\Gamma}: V^k(G)\longrightarrow V^{k-1/2}(\Gamma),\;
k\ge 0
\end{equation}
Evidently, these operators  are well-defined and bounded (see
\cite{T}).

\begin{deff}
{\rm A control $u(t,x)$ in (\ref{2.1})--(\ref{2.5}) is called
feedback \footnote{It will be clear later why defined control
really possesses feedback property} if
\begin{equation}\label{2.15}
v(t, \cdot)=\gamma_{\Omega}w(t, \cdot),\quad u(t,
\cdot)=\gamma_{\Gamma}w(t, \cdot) \qquad \forall t\ge 0
\end{equation}
where $(v(t, \cdot), u(t, \cdot))$ is the solution of
stabilization problem (\ref{2.1})--(\ref{2.5}) and $w(t, \cdot)$
is the solution of boundary value problem
(\ref{2.10})-(\ref{2.12}).}
\end{deff}

Evidently, if the solution $w$ of (\ref{2.10})-(\ref{2.12})
satisfies (\ref{2.39}), the pair $(v,u)$ defined in (\ref{2.15})
satisfies (\ref{2.6}). Since $(v,u)$ satisfies
(\ref{2.1})--(\ref{2.5}) as well, it forms a solution of the
initial stabilization problem (\ref{2.1})-(\ref{2.6}).


\subsection{Description of ``correct" initial conditions}

First of all we describe the set of initial conditions $\{ w_0\}$
such that solutions $w(t,x)$ of (\ref{2.10})-(\ref{2.12}) satisfy
(\ref{2.39}).

Let $G$ be domain (\ref{2.8}) and
\begin{equation}\label{2.16}
  V^0_0(G)=\{ v(x)\in V^0(G): (v,\nu )|_{\partial \Omega}=0\},\quad V^1_0(G)=V^1(G)\cap (H^1_0(G))^d
\end{equation}
where  $\nu (x)$ is the vector-field of outer unit normals to
$\partial G$. Evidently,
$$
    \|v\|_{V^0(G)}=\|v\|_{V^0_0(G)}:=\|v\|_{(L_2(G))^d};\quad \|v\|_{V^1_0(G)}:=\| \nabla v\|_{(L_2(G))^{d^2}}.
$$
Denote by
\begin{equation}\label{2.17}
  \hat \pi : (L_2(G))^2 \longrightarrow V^0_0(G)
\end{equation}
the operator of orthogonal projection. We consider the Oseen
steady state operator
\begin{equation}\label{2.18}
  Av\equiv -\hat \pi \Delta v+\hat \pi [(a(x),\nabla )v+(v,\nabla )a] :
  V^0_0(G) \longrightarrow V^0_0(G)
\end{equation}
and its adjoint operator $A^*$.  These operators are closed and
have the domain $  \cD (A)=V^2(G)\cap (H^1_0(G))^2$. Emphasize
that $\cD (A)$ consists of vector fields equal to zero on
$\partial G $. The spectrums $\Sigma (A), \Sigma (A^*)$ of
operators $A$ and $A^*$  are discrete subsets of a complex plane
$\mathbb{C}$
 which belong to a sector symmetric with respect to
$\bR $ and containing $\bR_+$. In other words, $A$ is a sectorial
operator.
 So spectrums $\Sigma (A), \Sigma (A^*)$
contain only eigenvalues of $A, A^*$, respectively. In virtue of
(\ref{2.9}) they are symmetric with respect to $\bR $, and
moreover $\Sigma (A)=\Sigma(A^*)$.

We rewritten the boundary value problem (\ref{2.10})-(\ref{2.12})
for Oseen equations   in the following form
\begin{equation}\label{2.33}
  \frac {dw(t,\cdot )}{dt}+Aw(t,\cdot )=0,\quad w|_{t=o}=w_0.
\end{equation}
where $A$ is the operator (\ref{2.18}). Then for each $w_0\in
V^0_0(G)$ the solution $w(t,\cdot)$ of {\rm (\ref{2.33})} is
defined by $w(t,\cdot)=e^{-At}w_0$ where $e^{-At}$ is the
resolving semigroup of problem (\ref{2.33}).

Let $\sigma >0$ satisfy:
\begin{equation}\label{2.35}
  \Sigma (A)\cap \{\lambda \in  \mathbb{C}:\Real\lambda=\sigma \}=
  \emptyset
\end{equation}
The case when there are certain points of $\Sigma (A)$ which are
in the left of the line $\{\Real\lambda=\sigma \}$ will be
interesting for us.

Denote by $X_\sigma^+(A)$ the subspace of $V^0_0(G)$ generated by
all eigenfunctions and associated functions of operator $A$
corresponding to all eigenvalues of $A$ placed in the set
  $\{ \lambda \in \mathbb{C}: \Real\lambda <\sigma\}$. By $X_\sigma^+(A^*)$
we denote analogous subspace corresponding to adjoint operator
$A^*$. We denote the orthogonal complement to $X_\sigma^+(A^*)$ in
$V^0_0(G)$ by $X_\sigma (A)\equiv X_\sigma $:
\begin{equation}\label{2.100}
   X_\sigma = V^0_0(G) \ominus  X_\sigma^+ (A^*)
\end{equation}
 One can show that subspaces $X_\sigma^+(A),\,   X_\sigma$ are invariant with
respect to the action of semigroup $e^{-At}$, and $X_\sigma
+X_\sigma^+(A)=V^0_0(G)$.

\begin{theorem}
Suppose that $A$ is operator {\rm (\ref{2.18})} and $\sigma >0 $
satisfies {\rm (\ref{2.35}).} Then for each $w_0 \in X_\sigma $
the  inequality (\ref{2.39}) holds. Besides, the solution of
problem (\ref{2.33}) with such initial conditions are defined by
the formula
\begin{equation}\label{2.34}
  w(t,\cdot )=e^{-At}w_0=(2\pi i)^{-1}\int_{\gamma}(A-\lambda I)^{-1}
  e^{-\lambda t}w_0 d\lambda .
\end{equation}
Here $\gamma$ is a contour belonging to $\rho
(A):=\mathbb{C}\setminus \Sigma (A)$ such that $\arg \lambda =\pm
\theta$ for $\lambda \in \gamma ,|\lambda |\ge N$ for certain
$\theta \in (0,\pi /2 )$ and for sufficiently large $N$.
Moreover, $\gamma$ encloses from the left   the part of the
spectrum $\Sigma (A)$ placed right of the line $\{\Real\lambda
=\sigma\}$. The complementary part of the spectrum $\Sigma (A)$ is
placed left of the contour $\gamma$.
\end{theorem}

Such contour $\gamma$ exists, of course.

Proof: See   \cite{F2}, \cite{F3}.


\subsection{Theorem on extension}

To complete the construction of stabilization for Oseen equations
(\ref{2.1}), (\ref{2.3}) we have to construct the operator $E$
extending initial condition $v_0$ from (\ref{2.3}) from $\Omega $
in $G$ such that $Ev_0:=w_0\in X_\sigma$.
 This $w_0(x)$ we take as initial value  in (\ref{2.12}). We consider here direct analog of construction from \cite{F6}.

 Introduce the space
\begin{equation}\label{4.19'}
V^0(\Omega ,\Gamma _0)=
       \{u(x)\in V^0(\Omega)\,:\, u\cdot \nu|_{ \Gamma_0}=0,\
       \exists\,v\in V^0_0(G):\  u=\gamma _{\Omega}v\}
 \end{equation}
supplied with the  norm:
$$
  \|u\|_{V^0(\Omega ,\Gamma _0)}=\inf_{w\in V^0_0(G):\gamma_\Omega w=u}\| w\|_{V^0_0(G)}
$$
where  $\gamma _\Omega $ is the restriction operator defined in
(\ref{2.91}) and
 $V^0_0(G)$ is defined in (\ref{2.16}).

We consider also the following closed subspace of $V^0_0(G)$:
$$
    V^0_0(G, \Omega )=\{ w\in V^0_0(G):\; w|_\Omega =0\}.
$$
For $u\in V^0(\Omega ,\Gamma_0)$ denote by $\widehat u$ the
element of the quotient space $V^0_0(G)/V^0_0(G, \Omega )$ such
that  for each $v\in \widehat u,\quad \gamma_\Omega v=u.$ Now we
define the extension operator
\begin{equation}\label{4.1010}
      L:\; V^0( \Omega ,\Gamma_0 )\to V^0_0(G)\quad \mbox{by}
\quad Lu\in \widehat u \quad \mbox{and}\quad
 Lu \bot_{V^0_0(G)}V^0_0(G, \Omega )
\end{equation}
(i.e. Lu is orthogonal to $V^0_0(G,\Omega )$ in the space
$V^0_0(G)$). Evidently, $\| L\| =1$.

 It is known (see \cite{F1}, \cite{F3}, \cite{F4}) that in the space $X^+_\sigma (A^*)$ one can choose a basis
$(d_1(x), \dots , d_K(x))$ such that restriction $(d_1(x)|_\omega,
\dots , d_K(x)|_\omega)$ on an arbitrary subdomain $\omega \subset
G$ forms a linear independent set of vector fields. We can define
space (\ref{2.100}) by the following equivalent form:
\begin{equation}\label{4.19}
   X_\sigma =\{v(x)\in V^0_0(G)\,:\,
  \int\limits_Gv(x)\cdot d_j(x)\,dx=0,\quad
  j=1,\dots, K\}.
\end{equation}

\begin{theorem}\label{t4.1}(\cite{F5}, \cite{F6})
There exists a linear bounded extension operator
\begin{equation}\label{4.20}
 E :V^0(\Omega ,\Gamma _0)\to X_\sigma,
\end{equation}
{\rm }i.e.\ $(Ev)(x)\equiv v(x)$ for $x\in\Omega $.
\end{theorem}

\begin{proof} Let subset $\omega_1\subset G\setminus \Omega$
be a domain with $C^\infty$- boundary $\partial \omega_1$ such
that $ \mbox{Int}(\partial \omega_1\cap \partial G)\ne \emptyset
$. In this set we consider the Stokes problem:
$$
    -\Delta w(x)+\nabla p(x)=v(x), \quad {\rm div}  w(x)=0,\quad x\in
 \omega_1 ;\quad w|_{\partial \omega_1}=0
 $$
As is well known, for each $v\in V^0(\omega_1)$ there exists a
unique solution $w\in V^1_0(\omega_1)\cap V^2(\omega_1)$ of this
problem. The resolving operator to this problem we denote as
follows: $(-\hat \pi \Delta)^{-1}_{\omega_1}v=w$. Extension of $(-\hat \pi
\Delta)^{-1}_{\omega_1}v$ from $\omega_1$ in $G$ by zero we also
denote as $(-\hat \pi \Delta)^{-1}_{\omega_1}v$. Evidently, $(-\hat \pi
\Delta)^{-1}_{\omega _1}v\in V^1_0(G)$.

We look for the extension operator $E$ in the form
\begin{equation}\label{4.23}
 E v(x)=(Lv)(x)+
            \bigg[\sum\limits^K_{j=1}c_j (-\hat \pi \Delta )^{-1}_{\omega _1} d_j\bigg](x),
\end{equation}
where $L$ is the operator (\ref{4.1010}),  $c_j$ are constants
which should be determined. Evidently, $E v(x)=v(x)$ if $x\in
\Omega$ for any $c_j$. Besides, $Ev \in V^1_0(G)$. To define
constants $c_j$ we note that by (\ref{4.19})  $E v\in X_\sigma$
if
\begin{equation}\label{4.25}
 \int\limits_Gd_k(x)
 \bigg[\sum\limits^K_{j=1}c_j(-\hat \pi \Delta )^{-1}_{\omega _1}d_j(x)\bigg]\,dx=
 -\int\limits_Gd_k(x)(Lv)(x))\,dx
\end{equation}
for $k=1,\dots,K$. As in \cite{F5}, \cite{F6} one can prove that
this system of linear equations has a unique solution.
\end{proof}

Aforementioned results imply the following result on stabilization
(see \cite{F1}-\cite{F5}).

\begin{theorem}\label{t4.6}
Let domains $\Omega $ and $G$ satisfy (\ref{2.8}). Then for each
initial value $v_0(x)\in V^1(\Omega,\Gamma _0)$ and for each
$\sigma >0$ there exists a feedback control~$u$ defined on~$\Sigma
$ such that the solution $v(t,x)$ of {\rm
(\ref{2.1})--(\ref{2.6})} satisfies the inequality
\begin{equation}\label{4.56}
   \|v(t,\cdot)\|_{V^0_0(\Omega)}\le ce^{-\sigma t}\ \ {\rm as}\ \ t\to\infty.
\end{equation}
\end{theorem}


\subsection{Real processes}\label{real}

Recall results from \cite{F6}, \cite{F7} on justification of
numerical simulation of stabilization construction described
above.  In virtue of definition (\ref{2.15})  stabilization
problem (\ref{2.1})-(\ref{2.6}) is reduced to problem
(\ref{2.10}), (\ref{2.12}), and (\ref{2.39}). Its numerical
simulation is what we have to justify.

 Let $e^{-At}=S(t)$ be resolving operator of
problem (\ref{2.10}), (\ref{2.12}).  Suppose that we calculate
this problem in   discrete time instants
$$
       t_1<t_2<\dots t_k<\dots
$$
where $t_k=k\tau$ and $\tau >0$ is fixed. Denote $S=S(\tau )$. Let
$\tilde w^k$ be the result of our calculations at time instant
$t_k$. Since numerical calculations can not be exact, we have
\begin{equation}\label{2.101}
            \tilde w^k=S\tilde w^{k-1}+\varphi^k
\end{equation}
where $\varphi^k$ is an error of calculation which is unknown for
us before time $t_k$. The sequence $\{ \tilde w^k\}$ defined by
(\ref{2.101}) is called real process.  We suppose that we can
estimate the error of our calculations a priori:
\begin{equation}\label{3.2}
 \| \varphi^{k} \|_{V^0_0(G)} \le \hat{\varepsilon } << 1,\; k> 0.
\end{equation}
where $\hat \varepsilon >0$ is a known quantity. Note also that at
time $t_k$ the vector $\tilde w^k$ is completely known (completely
observable) since it is result of our calculations.

Formulae (\ref{2.101}) is supplemented with initial condition
\begin{equation}\label{2.103}
             \tilde w^0=w_0,\qquad w_0\in X_\sigma .
\end{equation}
(We assume for simplicity that $\varphi^0=0$.)

Equations (\ref{2.101}), (\ref{2.103}) imply that
\begin{equation}\label{2.104}
             \tilde w^k=S^kw_0+\sum_{j=0}^{k-1}S^j\varphi^{k-j}
\end{equation}
In virtue of (\ref{2.104}), (\ref{2.103}) the estimate $\| \tilde
w^k\| \le ce^{-k\sigma \tau}$ is not true. Indeed, although
$S^kw_0\in X_\sigma$ and therefore $S^kw_0$ satisfies estimate of
such kind,  the fluctuations $\varphi^k$ possess nonzero
components belonging to $X_\sigma^+(A)$.  Hence, $\| S^j\varphi
^{k-j}\|_{V^0_0(G)}$ grows exponentially as $k\to \infty$ that
these terms destroy all stabilization construction described in
previous subsection.

To save this stabilization construction we use feedback mechanism
(see \cite{F6},\cite{F7}): in the time moment $t_k$ when $\tilde
w^k$ from (\ref{2.101})  is calculated we act on it by a special
projection operator $\Pi : V^0_0(G)\to X_\sigma$  that damps
undesirable   properties of fluctuations $\varphi^k$. This
operator has to satisfy the following properties:
$$
 i) \; \Pi \varphi =\varphi,\;  \forall \varphi \in X_\sigma , \quad
ii)  \; \forall \varphi \in V^0_0(G)\; (\Pi \varphi )(x)=\varphi
(x)\; \forall x\in \Omega .
$$
The operator satisfying these properties can be defined by the
formula analogous to (\ref{4.23}):
\begin{equation}\label{2.105}
 \Pi \varphi (x)=\varphi(x)+
            \bigg[\sum\limits^K_{j=1}c_j (-\pi \Delta )^{-1}_{\omega _1} d_j\bigg](x),
\end{equation}
where constants $c_j$ are defined by (\ref{4.25}) with $Lv=\varphi
$.

Applying to both parts of (\ref{2.101}), (\ref{2.103}) operator
$\Pi$. Taking into account that $X_\sigma$ is invariant with
respect to action of operator $S$ and using notation $w^k=\Pi
\tilde w^k$, we get the recurrent formulae for controlled real
process:
\begin{equation}\label{3.1}
            w^{k+1}=Sw^{k}+\Pi \varphi^{k+1}, \; \mbox{for}\;  k\ge 0;\qquad w^0=w_0
\end{equation}

In \cite{F6}, \cite{F7} the following estimate for controlled real
process has been proved:

\begin{theorem}\label{t4.7} 
Suppose that unpredictable fluctuations $\varphi^k$ satisfy
(\ref{3.2}), and operator $\Pi$ of projection on $X_\sigma$ is
defined by (\ref{2.105}).  Then controlled real process $\{ w^k\}$
defined by (\ref{3.1}) satisfies the estimate
\begin{equation}\label{4.220}
     \| w^k\|_{V^0(G)}\le (\gamma_0 e^{-\sigma k\tau }\| w_0\|_{V^0(G)}+
\| \Pi \| \hat \varepsilon /(1-e^{-\sigma \tau}))
\end{equation}
where  constant $\gamma_0\in (0,1)$ depends only on $\sigma $ if $\tau
>\tau_0$ and $\tau_0$ is a fixed magnitude \footnote{This property
of constant $\gamma_0$ is obtained in formula (\ref{4.10}) below.}, $\|
\Pi \|$ is the norm of operator (\ref{2.105}), and $\hat
\varepsilon >0$ is the magnitude from inequality {\rm
(\ref{3.2})}.
\end{theorem}

Note that outside a small neighborhood of the origin, estimate
(\ref{4.220}) is equivalent to estimate (\ref{2.39}). Since in
contrast to definition (\ref{2.34}) of solution $w(t,\cdot)$ to
problem (\ref{2.10}),(\ref{2.12}), in definition (\ref{3.1}) of
controlled real process fluctuations $\varphi^k$ permanently
arise, therefore estimates (\ref{4.220}), (\ref{2.39}) can not be
equivalent in a small neighborhood of the origin. Investigation of
behavior for $w^k$ in this neighborhood is the main goal of this
paper.





\section{Statistical problem in stabilization theory}


\subsection{Primary considerations}

We continue our investigation of the real process defined in
subsection \ref{real} above.

It is clear from the definition (\ref{2.101}), (\ref{2.103}) of
real process that its behavior is determined by its values in the
instants $t_k=k \tau$ when the unpredictable fluctuations
$\varphi^k $ arise. In the following we restrict ourselves by
considering the real process only in these points and thus we
obtain a discrete real process. So we consider iterated sequence
(\ref{3.1}) with  a fixed initial vector $w_0= w_0(x) \in
V_0^0(G)$ and unpredictable fluctuation $\varphi^{k+1} \in
V_0^0(G)$ satifying (\ref{3.2}).

It is   reasonable to assume that $ \varphi^{k }$ is a random
variable (for each $k$) defined on a probability space $(\Om,
\mathcal{A}, \bP)$, taking values in $V^0_0(G)$.  We also assume
that the random sequence $\varphi^{k }$ is independently
identically distributed (i.i.d.). Thus (\ref{3.1}) defines   a
random dynamical system (RDS).

For each $k$, $\phi^k$ has     distribution $\mu$, whose
probability measure $\mu (\om)$ is defined on the Borel
$\si-$algebra $\cB(V^0_0(G))$ of the space $V^0_0(G)$, and is
supported in a neighborhood of the origin:
\begin{equation}\label{3.3}
  \text{supp}\; \mu \subset B_{\hat{\eps}}
    =  \{ v\in V^0_0(G):\| v \|_{V^0_0(G)} \leq \hat{\eps} \}.
\end{equation}

Let us consider  the  random dynamical system (\ref{3.1}) from the
point of view of stabilization. The solution of the stabilization
problem can be derived in two stages: (a) To reach zero (or in
general setting, to reach a steady state); (b) To keep the
controlled solution $w^k$ near zero.

The solution of the stage (a) was explained above in Section 2. In
particular,  estimate (\ref{4.220}) of stabilized real process in
$G$ was obtained. Here and hereafter, we denote $\| w^k\| = \|
w^k\|_{V^0_0(G)}$, unless otherwise noted. The stage (a) is the part
of the process when $w^k$ does not reach the ball $B_{r_0}$ with
$r_0=\| \Pi\| \hat \varepsilon /(1-ce^{-\sigma \tau}) $ where $c\in
(0,1)$ is the constant from (\ref{4.220}).

Note that if $w^k \notin B_{r_0}$, then using (\ref{4.6}) (see below)
for estimating $S$ we get

 \begin{eqnarray}\label{3.5}
\|Sw^k - w^k\| > \|w^k\|-\|Sw^k\| \geq (1-\gamma_0e^{-\si \tau})\|w^k\|  \nonumber \\
     \geq  (1-\gamma_0e^{-\si \tau}) \frac{\|\Pi\| \hat{\varepsilon}}{1-\gamma_0e^{\si\tau}}
       = \|\Pi\| \hat{\varepsilon}.
\end{eqnarray}

This estimate means that at stage (a), the distance between $w^k$
and $Sw^k$ is more than the norm of the fluctuation $\Pi \phi^k$,
i.e., the deterministic component of random process $w^k$ is
prevailing. Therefore,  in stage (a) the behavior of $w^{k+1}$ is
determined mainly by the term $Sw^k$.

The situation in stage (b) when $w^k \in B_{r_0}$ is different.
Now the terms $Sw^k$ and $\Pi \phi^k$ from the right hand side of
(\ref{3.1}) have the equivalent order, and the motion of the
realization $w^k$ of the RDS (\ref{3.1}) becomes ``irregular or
chaotic".

The goal of this paper is to understand the behavior of the
controlled stabilized sequence $w^k$ in the stage (b). We will do
it with help of the modern theory of Markov chains or random
dynamical systems. We begin with some preliminaries.


\subsection{Gauss measures}

We recall some information about Gauss measures which will be used
below. Let $H$ be a Hilbert space with scalar product $(\cdot,
\cdot)$ and norm $\| \cdot \|$, and let $\mathcal{B}(H)$ be the
$\si-$algebra of Borel subsets of $H$. A measure $G(du)$ defined
on $\mathcal{B}(H)$ is called a Gauss measure if its Fourier
transform $\tilde{G}(v)$ is of the form
\begin{eqnarray}\label{3.6}
\tilde{G}(v) = \int e^{i(u,v)} G(du) = e^{i(v, a)-\frac12 (Kv,
v)}, \;\; v\in H,
\end{eqnarray}
where $a \in H$ is the mean vector, and $K: H\to H$ is a linear
self-adjoint positive \footnote{Positivness condition can be
weakened a little bit in the case of problem considered here: see
below footnote in the subsection \ref{s4.4}.}    trace class
operator, called correlation operator of $G$:
\begin{eqnarray}\label{3.7}
K^{\ast} = K > 0, \; \; \text{Trace}( K) = \sum_{j=1}^{\infty}
\lambda_j < \infty,
\end{eqnarray}
with $\{\lambda_j \}$ the set of eigenvalues of $K$.
Differentiation of (\ref{3.6}) with respect to $v$ implies that
$$
a = \int u G(du), \; (Kv_1, v_2)=(v_1,a)(v_2,a)-\int (v_1,
u)(v_2,u)G(du).
$$
Therefore $a$ is the mathematical expectation of the Gauss measure
$G$.

In particular, when $H=\bR^m$ and image $\Image K=\bR^m$, then for
each $\Gamma \in \mathcal{B}(\bR^m)$,
\begin{eqnarray}\label{3.8}
G(\Gamma) = \int p(dy) dy,
\end{eqnarray}
where density $p(y)$, $y\in \bR^m$, is defined by
\begin{eqnarray}\label{3.9}
p(y)=\frac1{(2\pi)^{m/2}\det K} \exp[-\frac12 (K^{-1}(y-a),
(y-a))].
\end{eqnarray}

Let $A:\; H\to H$ be a continuous map in Hilbert space $H$. As well-
known (\cite{VF}), each measure $\mu$ induces a new measure $A^{\ast}\mu (du)$
defined by
\begin{equation}\label{3.91}
A^{\ast}\mu (\om) = \mu (A^{-1} \om),\; \; \om \in \mathcal{B}(H),
\end{equation}
where $A^{-1}\om = \{x\in H: \; Ax\in \om \}$. This definition is
equivalent to
\begin{eqnarray}\label{3.10}
\int f(v)A^{\ast}\mu (dv) = \int f(Au) \mu(du)
\end{eqnarray}
for arbitrary $f$ for which at least one side in equation
(\ref{3.10}) is defined. Thus (\ref{3.6}) and (\ref{3.10}) imply
that if $\mu =G$ is the Gauss measure defined above
and map $A$ is linear, then
$A^{\ast}G$ is a Gauss measure with mathematical expectation $a_1
= Aa$ and correlation operator $K_1=AKA^{\ast}$.


\subsection{ Distribution of $\varphi^k$}

We consider the right hand side of (\ref{3.1}) where $\phi^k$ is
an i.i.d. random sequence. We suppose that the distribution
$\mathcal{D}(\phi^k)$ of $\phi^k$ has the form
\begin{eqnarray}\label{3.11}
\mathcal{D}(\phi^k) = c \chi_{\hat{\eps}}(u)G(du) = \nu (du),
\end{eqnarray}
where $c=(\int_{B_{\hat{\eps}}}G(du))^{-1}$, and
\begin{equation}\label{3.12}
\chi_{\hat{\eps}}(u)=
    \begin{cases}
    1, &\|u\| \leq \hat{\eps} \\
    0, &\|u\| > \hat{\eps}
    \end{cases}
\end{equation}
is the characteristic function of the ball $B_{\hat{\eps}}$
defined in (\ref{3.3}), $G(du)$ is the Gauss measure with
mathematical expectation $a=0$ and correlation operator $K$
satisfying  conditions (\ref{3.7})

In virtue of definition (\ref{3.11}),
\begin{eqnarray}\label{3.14}
\nu (\om) = cG(B_{\hat{\eps}}\cap \om),\;\; \om \in
\mathcal{B}(V_0^0(G))
\end{eqnarray}
and therefore $\nu(\om)$ is a probability measure on
$\mathcal{B}(V_0^0(G))$ supported on the ball $B_{\hat{\eps}}$.


\subsection{The main result}

Since $\phi^k$ are i.i.d., $\Pi\phi^k$ are i.i.d. as well, and
$D(\Pi\phi^k) \in X_{\si}$. Thus RDS (\ref{3.1}) defines a family
of Markov chains in $X_{\si}$ with transition function
\begin{eqnarray}\label{3.15}
P(k, w_0, \Gamma)= \bP \{w^k(w_0) \in \Gamma \},\;\;\; \Gamma \in
\cB(X_{\si}),
\end{eqnarray}
where $w^k = w^k(w_0)$ is defined by (\ref{3.1}) and $\bP$ is
probability measure defined on $\sigma$-algebra $\mathcal{A}$ of
subsets of probability space $\Omega$. Let $\cP(X_{\si})$ be the
space of Borel probability measures on $\cB(X_{\si})$ and
$C_b(X_{\si})$ be the space of continuous bounded functions on
$X_{\si}$. Moreover, we denote by $\mathfrak{P}_k$ and
$\mathfrak{P}_k^{\ast}$ the corresponding Markov semigroups acting
in $C_b(X_{\si})$ and $\cP(X_{\si})$, respectively:
$$
\fP_k f(v)=\bE f(w^k(w_0))\equiv \int f(z)P(k, w_0, dz),\;\; f\in
C_b(X_{\si}),
$$
$$
\fP_k^{\ast} \mu(\Gamma) = \int_H P(k, w_0, \Gamma)\mu(dw_0),\;\;\;
\mu\in \cP(X_{\si}),
$$
where $\bE$ is for the mathematical expectation, and $P(k, w_0,
\Gamma)$ is defined in (\ref{3.15}).

A continuous function $f(u)$ on $X_{\si}$ is called Lipschitz if
$$
\sup_{v\in X_{\si}, \|v\| \leq 1} \frac{|f(u+v)-f(u)|}{\|v\|}
\equiv \text{lip} f(u) < \infty, \; u\in X_{\si}.
$$
We denote $\text{Lip} f = \|\text{lip} f(\cdot)\|_{C(X_{\si})}$.

A measure $\mu \in \cP(X_{\si})$ is called a stationary measure
for the RDS (\ref{3.1}) if $\fP_k^{\ast} \mu = \mu, \forall k$.
The main theorem of this paper is  as follows.

\begin{theorem} \label{main}
The random dynamical system (\ref{3.1}) has a unique stationary
measure $\hat{\mu}$. Moreover, there exists a constant $\gamma \in
(0, 1)$ such that
\begin{eqnarray}\label{3.16}
|\int f(z) P(k, w_0, dz)-\int f(z)\hat{\mu}(dz)| \leq c\gamma^k,
\; k=1, 2, 3, \cdots,
\end{eqnarray}
for every Lipschitz function $f$ on $X_{\si}$ such that
$\|f\|_{C_b(X_{\si})} \leq 1$ and $\text{Lip} f \leq 1$. The
constant $c$ depends only on initial state $\|w_0\|$.
\end{theorem}

Note that the uniqueness of the stationary measure means that the
RDS (\ref{3.1}) is ergodic. The exponential convergence
(\ref{3.16}) means that RDS (\ref{3.1}) possesses the property of
exponential mixing.

Theorem \ref{main} provides us the possibility of calculating
easily the probability characteristics of Markov chain
(\ref{3.1}). Indeed, in numerical simulation we actually obtain
certain realizations $w^k=w^k(\om)$ of RDS (\ref{3.1}), where $\om
\in \Om$ is a random sample. By the strong law of large numbers
\begin{eqnarray}\label{3.17}
\lim_{N\to\infty} \frac1{N}\sum_{k=0}^{N} w^k(\om) \to \int
w\mu(dw), \; \text{as} \;  N\to \infty,
\end{eqnarray}
where $\mu(dw)$ is an invariant measure of RDS (\ref{3.1}). In
virtue of Theorem 3.1, $\mu(dw)=\hat{\mu}(dw)$. So while
calculating $w^k(\om)$, we can simultaneously obtain mathematical
expectation of $\hat{\mu}$. Moreover, (\ref{3.16}) gives us the
convergence rate in (\ref{3.17}).

The topic connected with (\ref{3.17}) will be studied in detail in
some other place. Here we only note that strong law of large
numbers was derived from ergodicity of Navier-Stokes equations in
\cite{K2} in the case when random force is white noise. Of course
in the case of kick forces this derivation can be made as well.


\subsection{Ergodic theorem} \label{ergodic}

In order to prove Theorem 3.1, we use a result   in \cite{KSh, K}
on ergodicity of a Markov chain, proved by coupling techniques.
Let us formulate this result in a form suitable for our problem.
Consider a Markov chain (or a RDS) in a Hilbert space $H$
\begin{eqnarray}\label{3.18}
u(k) = T u(k-1) + \eta_k, \; u(0)=u_0,
\end{eqnarray}
where $u_0 \in H$, $T: H\to H$ is a linear bounded operator
  such that
\begin{eqnarray}\label{3.19}
 \|Tu\| \leq \gamma_0 \|u\|, \; \forall u\in H,
\end{eqnarray}
for a constant $\gamma_0\in (0,1)$.

Assume also that there exists an orthonormal  basis $\{ e_j\}$ in
$H$ (true for any separable Hilbert space $H$) and a sequence of
subspaces $H_1\subset H_2\subset \dots H_k \subset \dots$ such
that
$$
     H_k=\mbox{span}\{ e_1, \dots ,e_{r_k}\} \quad \mbox{where}\quad r_k\rightarrow \infty \quad \mbox{as}
\quad k\rightarrow \infty
$$
Denote by $H^\perp_k$ the orthogonal complement to $H_k$ in $H:\;
H=H_k\oplus H_k^\perp$, and designate by
$$
        Q_k : H\longrightarrow H_k, \qquad Q_k^\perp : H\longrightarrow H^\perp_k
$$
orthogonal projectors.  Suppose that
\begin{eqnarray}\label{3.20}
 \|Q_k^\perp Tu\| \leq \gamma_k \|u\|, \; \forall u\in H, \quad \mbox{and}\; \gamma_k \rightarrow 0 \; \mbox{as}
\; k\rightarrow \infty
\end{eqnarray}

At last  assume that in (\ref{3.18}), $\eta_k$ is an i.i.d. random
sequence with distribution $\mu(\om), \om\in H$, such that the
projection $Q_k^* \mu$ on $H_k$, has a continuous density $R(x)$
with a compact support:
$$
Q^*_k \mu (dx) = R(x)dx
$$
where $H_k \ni x=\sum x_i e_i, \; dx=dx_1\cdots dx_{r_k}$.
Moreover, this density $R(x)$ satisfies the condition:
 \begin{eqnarray}\label{3.201}
  \int_{H_k}|R(x-v_1)-R(x-v_2)|\, dx\le c\| v_1-v_2\|_{H_N}
\end{eqnarray}
where the constant $c>0$ does not depend on $v_1,v_2\in H_N$.

Under aforementioned assumptions in \cite{KSh, K}, the following
theorem has been proved (see \cite{K}).

\begin{theorem} \label{kuksin}
The RDS (\ref{3.18}) has a unique stationary measure $\hat{\mu}$.
Moreover, there exists a constant $\gamma \in (0, 1)$ such that
$$
|\int f(z) P(k, u, dz)-\int f(z)\hat{\mu}(dz)| \leq c\gamma^k,\;
k=1,2,\cdots,
$$
for every Lipschitz function $f$ on $H$ satisfying
$|f|_{C_b(H)}\leq 1$ and $ \text{Lip} f \leq 1$. Here $P(k,u,dz)$
is the transition function (see (\ref{3.15})) corresponding to RDS
(\ref{3.18}).
\end{theorem}


\section{Check of assumptions (\ref{3.19}), (\ref{3.20})}

To prove Theorem \ref{main}, we have to check that RDS (\ref{3.1})
satisfies
 conditions (\ref{3.19}), (\ref{3.20}), and (\ref{3.201}) of Theorem \ref{kuksin}.
In this section we check the first and the second of them.


\subsection{Subspaces of $V_0^0(G)$}\label{s4.1}

 First we introduce an orthogonal decomposition of
$V_0^0(G)$ in order to define analogs of subspaces $H_k$ from
subsection \ref{ergodic}.  Recall that an important subspace  of
$V_0^0(G)$ is (see definition in (\ref{4.19}))

\begin{eqnarray}\label{4.1}
  X_{\si}=\{v\in V_0^0(G): \int_G v(x)d_j(x) dx=0, j=1,\cdots,m\},
\end{eqnarray}
where $\{d_j(x), j=1,\cdots,m\}$  is the basis of
$X^+_{\si}(A^{\ast})$, constructed in \cite{F2} from eigenvectors
and associated vectors
 of the operator $A^{\ast}$,
corresponding to all eigenvalues $\lambda_j$ satisfying $Re
\lambda_j > \si$.

Let $\si < \si_1<\dots <\si_k \rightarrow \infty$ as $k\to \infty$
be  a sequence of numbers  which satisfy, as $\si$, the condition
\begin{eqnarray}\label{4.2}
  \{\lambda \in \bC: Re \lambda = \si_k\}\cap \sigma (A) =\emptyset \qquad \forall \; k.
\end{eqnarray}
Analogously to (\ref{4.1}), we can define the spaces
\begin{eqnarray}\label{4.3}
  X_{\si_k}=\{v\in V_0^0(G): \int_G v(x)d_j(x) dx=0, j=1,\cdots,n_k\},
\end{eqnarray}
where basis $\{ d_j(x), j=1,\dots ,n_k\}$ in $X^+_{\sigma_k}(A^*)$
is constructed from eigenvectors and associated vectors
 of the operator $A^{\ast}$,
corresponding to all eigenvalues $\lambda_j$ satisfying $Re
\lambda_j > \si_k$. This basis is an extension of the basis in
$X^+_\sigma (A^*)$ from (\ref{4.1}), and it is constructed by the
same rules as the basis from (\ref{4.1}).

Since $\sigma_i>\sigma_j >\si$ for $i>j$ we have that  $n_i>n_j>m$
and therefore $X_{\si_i} \subset X_{\si_j}\subset X_\si$.
Moreover, we introduce  the following subspaces of $V_0^0(G)$. Let
$X_{\si\si_1}$ be an orthogonal complement in $X_{\si}$ for the
subspace $X_{\si_1}$, and $X_{\si_k\si_{k+1}}$ be an orthogonal
complement in $X_{\si_k}$ for the subspace $X_{\si_{k+1}}$. In
other words, $X_{\si\si_1}$ and $X_{\si_k\si_{k+1}}$   are
subspaces satisfying
\begin{eqnarray}\label{4.4}
  X_{\si_1} \oplus  X_{\si\si_1}= X_{\si}\qquad X_{\si_{k+1}} \oplus  X_{\si_k\si_{k+1}}= X_{\si_k}.
\end{eqnarray}
We also define the subspace $X_{\si}^{\perp} \subset V_0^0(G)$,
which is the orthogonal complement of $X_{\si}$ in $V_0^0(G)$:
\begin{eqnarray}\label{4.5}
  X_{\si} \oplus  X_{\si}^{\perp}
=  V_0^0(G).
\end{eqnarray}
Evidently, $X^\perp_\si=X^+_\si(A^*)$. The subspace $ X_{\si} $
will play the role of space $H$ in subsection \ref{ergodic}.
Likewise,
\begin{eqnarray}\label{4.51}
 X_{\si \si_k}:=X_{\si\si_1}\oplus \dots \oplus X_{\si_{k-1}\si_k}
\end{eqnarray}
  will play the role of $H_k$ and $X_{\si_k}$ will play the role of $H_k^{\perp}$.
Recall that subspace  $X_{\si} $  is invariant with respect to the
operator $S$ in RDS $(\ref{3.1})$. That is why we put in
(\ref{3.18})
$$
T=S|_{X_{\si}}.
$$
Now we construct the basis $\{e_j\}$ from subsection
\ref{ergodic}. Let $\{e_1, \dots ,e_m\}$ be orthogonalisation of
basis  $d_1,\dots ,d_m$ in $X^+_\si(A^*)=X^\perp_\si$. Evidently,
$\{e_1, \dots ,e_m\}$ forms a  orthonormal basis in $X^\perp_\si$.
Continuing orthogonalization process for
$$
\{ d_{m+1},\dots ,d_{n_1}\} ,  \dots , \{ d_{n_{k-1}+1},\dots
,d_{n_k}\} ,  \dots
$$
 we get orthnormal basis
$\{ e_{m+1},\dots ,e_{n_1}\}$ in $X_{\si \si_1}$ and $\{
e_{n_{k-1}+1},\dots ,e_{n_k}\}$ in $X_{\si_{k-1} \si_k}$ for
$k=2,3,\dots$. We have to prove now that countable orthonormal
system $\{ e_1,\dots , e_j, \dots \}$ forms a basis in $V^0_0(G)$.
For this it is enough to establish that this system is dense in
$V^0_0(G)$. But in virtue of Keldysh Theorem (see \cite{Kel},
\cite{GK}), the system $\{ d_j, j\in \mathbb{N}\}$ constructed by
eigenfunctions and associated functions of operator $A^*$ adjoint
to operator (\ref{2.18}) is dense in $V^0_0(G)$. Hence, system $\{
e_j, j\in \mathbb{N}\}$ obtained from
 $\{ d_j, j\in \mathbb{N}\}$ by orthogonalization process is also dense.
 Therefore the  system
$\{ e_j, j=m+1, \dots , m+k, \dots \}$ forms basis in the space
$H=X_\si$.

In the next  subsection we establish inequality (\ref{3.19}).


\subsection{Certain properties of RDS (\ref{3.1}) }

We prove the following assertion.
\begin{lem}\label{lem4.10}
Let $S$ be the operator in (\ref{3.1})and $\si >0$ be given. Then
for each $\gamma_0 \in (0,1)$, there exists $\tau > 0$ such that
for $S=S(\tau)$ the following estimate holds
\begin{eqnarray}\label{4.6}
  \|Su\| \leq \gamma_0 \|u\|,\; \forall u\in X_{\si}.
\end{eqnarray}
\end{lem}

\begin{proof}
Recall that the basic space is $V_0^0(G)$ and therefore we use the
notation $\|\cdot\| = \|\cdot \|_{V_0^0}$. Inequality
(\ref{4.6})  follows from the bound established in
\cite{F2}:
\begin{eqnarray}\label{4.62}
 \| \int_{\gamma_{\si}} (A-\lambda I )^{-1} e^{-\lambda \tau} d\lambda \| =
\| \int_{-\gamma_{\si}} (A+\lambda I )^{-1} e^{\lambda \tau}
d\lambda \|
   \leq c e^{-\si \tau},
\end{eqnarray}
where $A$ is the infinitesimal generator for $S(t)$, and the
contour $-\gamma_{\si}$ is defined as follows:
$$
     -\gamma_{\si}=\gamma_{\si}^1\cup \gamma_{\si}^2
$$
where
$$
\gamma_{\si}^1=\{\lambda \in \bC, Re \lambda = -\si, Im \lambda
\in [-(\si + \theta)\tan(\pi-\psi),(\si + \theta)\tan(\pi-\psi) ]
\}
$$
\begin{eqnarray}\label{4.91}
\gamma_{\si}^2=\{\lambda \in \bC, Re \lambda <-\si, \lambda=\gamma
e^{\pm i\psi} + \theta, \text{for}\; \gamma\in [\frac{\si
+\theta}{|\cos\psi |}, \infty)\}
\end{eqnarray}
with $\theta >0$ and $\pi/2 <\psi <\pi$   fixed. We have to prove
that $c$ in (\ref{4.62}) can be chosen independent of $\si
>0$. Since $\gamma_{\si}$ belongs
to the resolvent set of the operator $A$, we can get as in
\cite[Lemma 4.7]{F2} that
\begin{eqnarray}\label{4.9}
  \|(\lambda I + A)^{-1} \| \leq \frac{M_1}{1+|\lambda|},
       \qquad  \lambda \in -\gamma_{\si},
\end{eqnarray}
with $M_1>0$ independent of $\lambda \in -\gamma_\sigma$. We see
that
 $$
\| \int_{-\gamma_{\si}} (\lambda I + A)^{-1} e^{\lambda
\tau}d\lambda \| \le I_1+I_2
$$
where
\begin{equation}\label{4.92}
   I_1= \int_{\gamma_{\si}^1} \| (\lambda I + A)^{-1}\| \;
  | e^{\lambda \tau}d\lambda | \quad
  I_2=\int_{\gamma_{\si}^2} \| (\lambda I + A)^{-1}\| \;
  | e^{\lambda \tau}d\lambda |
\end{equation}
Using (\ref{4.9}) we get
\begin{align*}
I_1 &\le \int_{-(\si + \theta)\tan(\pi-\psi)}^{(\si +
\theta)\tan(\pi-\psi)} \frac{M_1
e^{-\si\tau}dx}{1+\sqrt{\si^2+x^2}} \leq M_1
e^{-\si\tau}\int_{-(\si + \theta)\tan(\pi-\psi)}^{(\si +
\theta)\tan(\pi-\psi)} \frac{d x/\si}{\sqrt{1+(x/\si)^2}} \\
& =M_1e^{-\si\tau}\int_{-(1 + \theta/\si)\tan(\pi-\psi)}^{(1 +
\theta/\si)\tan(\pi-\psi)} \frac{dy}{\sqrt{1+y^2}}
 \leq M_1ce^{-\si\tau},
\end{align*}
where $c$ does not depend on $\si \geq 1$ and $\tau >0$. If $\si
\in (0, 1)$, then
$$
I_1 \leq  \int_{-(1 + \theta)\tan(\pi-\psi)}^{(1 +
\theta)\tan(\pi-\psi)} \frac{M_1e^{-\si\tau}dx}{1+|x|} =
M_1ce^{-\si\tau},
$$
with $c$ not depending on $\si \in (0,1)$ and $\tau >0$. Moreover,
by change of variables $x=(\gamma |cos\psi|-\theta)$, we get
\begin{align*}
I_2 &\leq 2M_1  \int_{\frac{\si + \theta}{|\cos\psi|}}^{\infty}
\frac{\exp[-(\gamma|\cos\psi|-\theta)\tau]
d\gamma}{1+\sqrt{(\gamma|\cos\psi|-\theta)^2+\gamma^2\sin^2\psi}}
\\
&=\frac{2M_1}{|\cos\psi|}\int_{\si}^{\infty} \frac{\exp[-x\tau]
dx}{1+\sqrt{x^2+ (x + \theta)^2\tan^2\psi}} \leq
\frac{2M_1}{|\cos\psi|}\int_{\si}^{\infty} \frac{\exp(-x\tau)
dx}{1+x}
\\
&= \frac{2M_1}{|\cos\psi|}\int_{\si\tau}^{\infty} \frac{\exp(-y)
dy}{\tau + y}
=\frac{2M_1\exp(-\si\tau)}{|\cos\psi|}\int_0^{\infty}
\frac{\exp(-z) dz}{\tau(1+\si) + z},
\end{align*}
where in the final step, we used the change of variables
$y=z+\si\tau$.

Note that
$$
       \int_0^\infty\frac{\exp(-z)\,dz}{\tau (1+\si)+z}\le 1 \qquad \mbox{if} \quad \tau \ge 1
$$
For $\tau \in (0,1)$ we have
$$
\int_0^{\infty} \frac{\exp(-z) dz}{\tau(1+\si) + z} \le
\int_0^{\infty} \frac{\exp(-z) dz}{\tau+z} \le
\int_0^{1-\tau}\frac{dz}{\tau +z} +\int_0^\infty e^{-z} dz=1-\ln
\tau
$$
So we have
\begin{eqnarray}\label{4.93}
I_2 \leq \frac{2M_1}{|\cos\psi|} e^{-\si\tau} c_1,
\end{eqnarray}
where $c_1 $ does not depend on $\tau > \tau_0$ and $\si >0$.

Thus we have proved that $c$ in (\ref{4.62}) does not depend on
$\si >0$ and $\tau >\tau_0$. So for $u\in X_{\si}$,
\begin{eqnarray} \label{4.10}
\| S(\tau)u\| \leq cM_1 e^{-\si \tau}.
\end{eqnarray}
where $c>0$ does not depend on $\sigma >0$ and $\tau >\tau_0$.
Therefore for given $\gamma_0$ and $ \si>0$, we can take $\tau>0$
such that $c e^{-\si\tau}=\gamma_0$.
\end{proof}

Inequality (\ref{3.19}) evidently follows from (\ref{4.6}).


\subsection{Estimate $S(\tau)$ on $X_{\si_k}$ for large $k$}

We check here bound (\ref{3.20}) for $T=S(\tau)$ and
$Q^\perp_kH=X_{\si_k}$. Simultaneously we make more precise the
choice of $\si_k$. Together with operator $A$ defined in
(\ref{2.18}) we consider the Stokes operator
\begin{eqnarray}\label{4.120}
  A_0=-\hat \pi \Delta : V^0_0(G)\rightarrow V^0_0(G)
\end{eqnarray}
where $\hat \pi$ is projector (\ref{2.17}). Operator $A_0$ is
positive self-adjoint with domain
$\mathcal{D}(A_0)=\mathcal{D}(A)=V^2(G)\cap V^1_0(G)$ and with
discrete spectrum. Let $\{ \varepsilon_j\}$ be eigenvectors of
$A_0$ forming an orthonormal basis in $V^0_0(G)$ and $0<\mu_1 \le
\mu_2 \le \dots $ be corresponding eigenvalues,     taking into
account of their multiplicities. Then for each $q\in \mathbb{R}$
the power $A^q_0$ can be defined by the formula
$A^q_0v=\sum\mu_j^q(v,\varepsilon_j)_{V^0_0(G)}\varepsilon_j $.
Eigenvalues $\mu_j$ possess the following asymptotic for large
$j$:
\begin{eqnarray}\label{4.121}
   \mu_j=\beta_0j^{2/d}+O(j^{2/d}/\ln{j})\qquad \mbox{as}\quad j\rightarrow \infty
\end{eqnarray}
where $\beta_0>0, d$ is dimension of $G$ (i.d. $d=2$ or $3$),
$O(j^{2/d}/\ln{j})$ is a function satisfying the estimate
$|O(j^{2/d}/\ln{j})|\le cj^{2/d}/\ln{j}$ as $j\to \infty$ with
$c>0$ independent on $j$.  Asymptotics (\ref{4.121}) was obtained
by K.I.Babenko \cite{B}. \footnote{Actually, K.I.Babenko proved (\ref{4.121}) in
\cite{B} in the case d=\mbox{dim}G=3 only. But his proof can be extended in the
case d=2 as well.}

We can write operator $A$ from (\ref{2.18}) as follows:
\begin{eqnarray}\label{4.122}
  A=A_0+A_1
\end{eqnarray}
 where $A_1v=\hat \pi [(a(x),\nabla)v+(v,\nabla)a]$. Using \cite[Ch.3, Lemma 4.5]{F8} one
can easily get the following bound:
\begin{eqnarray}\label{4.123}
  \|A_1A_0^{-1/2}\|=b<\infty
\end{eqnarray}

Now we describe one result of M.S.Agranovich announced in
\cite[Bound (6.61)]{A1}. Although it is obtained for general
abstract operators we formulate it in the case of Oseen operator.
Its proof will be published in \cite{A2}.This result consists of
the choice of sequence $\si_k \to \infty $ as $k\to \infty$ such
that on segments
\begin{eqnarray}\label{4.124}
         \Gamma_k=\{ \lambda =\si_k+i\gamma ;\; |\gamma|\le b'\si_k^{1/2}\},\quad b'>0
\; \mbox{does not depend on}\; k,
\end{eqnarray}
resolvent $(A-\lambda I)^{-1}$ possesses some optimal bound.
Numbers $\si_k$ are found on segments \footnote{Actually,
Agronovitch's result was proved in \cite{A2} under assumption that
reminder term in (\ref{4.121}) has the form $O(j^r)$ with $r<2/d$
(that is stronger than $O(j^{2/d}/\ln{j})$ from (\ref{4.121})),
and $\si_k$ are looked for in segments
$[k^{2\rho/d},(k+1)^{2\rho/d}]$ with some $\rho>4$. But if we
change these segments on segments (\ref{4.125}) then a
streightforward repeating of the corresponding proof from
\cite{A2} leads to the desired estimate.}
\begin{eqnarray}\label{4.125}
          \Delta_k=[e^{2k/d},e^{2(k+1)/d}], \qquad d=\mbox{dim}G=2\; \mbox{or} \; 3.
\end{eqnarray}

\begin{lem}\label{lemAgr} Suppose that an operator $A$ has the form
(\ref{4.122}) where $A_0$ is self-adjoint positive operator with
discrete spectrum and eigenvalues satisfying (\ref{4.121}), and
$A_1$ satisfies (\ref{4.123}). Then there exists $\si_k\in
\Delta_k$ such that  for $\lambda \in \Gamma_k$ the following
estimate holds:
\begin{eqnarray}\label{4.126}
   \| (A-\lambda I)^{-1}\| \le c_1c_2^{\si_k^{\frac{d-1}{2}}}\si_k^{-1/2}
\end{eqnarray}
where segments $\Gamma_k, \Delta_k$ are defined in
(\ref{4.124}),(\ref{4.125}), respectively.
\end{lem}

Using Lemma \ref{lemAgr} we can check (\ref{3.20}). Recall that
$S(t)=e^{-At}$ is the resolving semigroup of problem (\ref{2.33})
where $A$ is Oseen operator (\ref{2.18}), and $S=S(\tau)$ where
$\tau$ is a fixed number chosen in Lemma \ref{lem4.10} such that
(\ref{4.6}) is true. Note that if we would increase $\tau$,
(\ref{4.6}) is true as well. Recall that the space $X_{\si_k}$
defined in (\ref{4.3}) is invariant with respect to the operator
$S=S(\tau)$.

\begin{theorem}
Let $A$ be Oseen operator (\ref{2.18}), $S=S(\tau)=e^{-A\tau}$,
and sequence $\si_k \to \infty$ as $k\to \infty $  be chosen in
Lemma \ref{lemAgr}. Then there exists $\tau_0>0$ such that for
each $\tau >\tau_0$ on the spaces $X_{\si_k}$ the following
estimate hold:
\begin{eqnarray}\label{4.127}
  \| S(\tau)u\|_{X_{\si_k}}\le \gamma_k \| u\|_{X_{\si_k}}, \qquad \mbox{where}\quad
\gamma_k \rightarrow 0\quad \mbox{as}\quad k\rightarrow \infty,
\end{eqnarray}
where $\gamma_k's$ do not depend on $u\in X_{\si_k}$.
\end{theorem}

\begin{proof}
It is clear that for $u\in X_{\si_k}$
\begin{eqnarray*}
  S(\tau)u=(2\pi i)^{-1} \int_{\gamma_{\si_k}} (A-\lambda I )^{-1}u e^{-\lambda \tau} d\lambda  =
 -(2\pi i)^{-1}\int_{-\gamma_{\si_k}} (A+\lambda I )^{-1}u e^{\lambda \tau} d\lambda
\end{eqnarray*}
where $-\gamma_{\si_k}=\gamma^1_{\si_k}\cup \gamma^2_{\si_k}$ and
similarly to (\ref{4.91})
\begin{eqnarray*}
\gamma_{\si_k}^1=\{\lambda \in \bC, Re \lambda = -\si_k, Im
\lambda \in [-(\si_k +
\theta)\tan(\pi-\psi),(\si_k + \theta)\tan(\pi-\psi) ] \} \\
\gamma_{\si_k}^2=\{\lambda \in \bC, Re \lambda <-\si_k,
\lambda=\gamma e^{\pm i\psi} + \theta, \text{for}\; \gamma\in
[\frac{\si_k +\theta}{|\cos\psi |}, \infty)\}
\end{eqnarray*}
Doing calculation as in proof of Lemma \ref{lem4.10} we get the
same formulas where $\si$ is changed to $\si_k$ only. Then
estimation of the term $I_2$ gives as in (\ref{4.93}):
\begin{eqnarray}\label{4.128}
I_2 \leq \frac{2M_1}{|\cos\psi|} e^{-\si_k\tau} c_1,
\end{eqnarray}
where $c_1=c_1(\tau_0)$ does not depend on $\tau >\tau_0$. To
estimate $I_1$ we use (\ref{4.126}) instead of (\ref{4.9}). More
precisely we use the estimate
\begin{eqnarray}\label{4.129}
   \| (A+\lambda I)^{-1}\| \le c_1e^{\si_k \ln{c_2}}\si_k^{-1/2},\quad
\lambda \in -\tilde \Gamma_k
\end{eqnarray}
where $\tilde \Gamma_k=\{ \lambda =\si_k+i\gamma ;\; |\gamma|\le
(\si_k+\theta) \tan (\pi -\psi)\}$. If $\lambda \in -\Gamma_k$,
estimate (\ref{4.129}) directly follows
 from (\ref{4.126}). For $\lambda \in -\{ \tilde \Gamma_k \setminus \Gamma_k\}$ situation is easier and one can prove estimate similar to (\ref{4.9}) in right side of which $|\lambda |$
is changed on $|\lambda |^{1/2}$ that is stronger than
(\ref{4.129}). (This has been done in \cite{A2}).  Applying
(\ref{4.129}) to $I_2$ defined in (\ref{4.92}) (with $\si$ changed
on $\si_k$) we get:
\begin{eqnarray}\label{4.130}
   I_1\le \int_{-(\si_k+\theta)\tan{(\pi -\psi)}}^{(\si_k+\theta)\tan{(\pi -\psi)}}
c_1e^{-\si_k(\tau -\ln{c_2})}\si^{1/2}\,dx \le \tilde c
e^{\si_k(\tau-\ln{c_2})}
\end{eqnarray}
where $\tilde c$ does not depend on $k$. If we choose $\tau
>\ln{c_2}$ then (\ref{4.128}), (\ref{4.130}) imply that
$$
     I_1+I_2\le \gamma_k \rightarrow 0\qquad \mbox{as} \quad k\rightarrow \infty .
$$
This proves (\ref{4.127}).
\end{proof}

Inequality (\ref{3.20}) follows from (\ref{4.127}).


\section{Reduction to the finite dimensional case}

Now we need only to check the condition (\ref{3.201}). The rest
part of the paper is devoted to prove that the distribution
$\Pi\phi^k$ of random forces in  RDS (\ref{3.1}) satisfies this
property. This will then complete this paper. First, we project
$\Pi\phi^k$ on finite-dimensional subspace.

\subsection{Calculation of the projection for probability
distribution $\Pi\phi^k$}

In this subsection we calculate probability distribution for
random variable $Sw^k+\Pi\phi^{k+1}$ in RDS (\ref{3.1}) under
assumption that $Sw^k$ is a fixed vector. Since $Sw^k\in X_{\si}$
and $\Pi: V_0^0\to X_{\si}$ is a projection on $X_{\si}$, the
probability distribution $\mathcal{D}(Sw^k+\Pi\phi^{k+1})$ is
supported on $X_{\si}$. It is enough to calculate
$\mathcal{D}(\Pi\phi^{k+1})$ because
$\mathcal{D}(Sw^k+\Pi\phi^{k+1})$ is simply the shift of
$\mathcal{D}(\Pi\phi^{k+1})$ along the vector $Sw^k$. In virtue of
(\ref{3.11})-(\ref{3.14}), $\mathcal{D}(\phi^{k+1})=\nu$, where
the measure $\nu$ is defined by
\begin{eqnarray} \label{4.11}
\nu(\om)=cG(B_{\hat{\eps}}\cap \om),\; \forall \om \in \cB(V_0^0),
\end{eqnarray}
where $G$ is the Gauss measure with mathematical expectation $a=0$
and correlation operator $K$ satisfying (\ref{3.7}). Clearly,
$\mathcal{D}(\Pi\phi^{k+1}) = \Pi^{\ast}\nu$. Since
\begin{eqnarray} \label{4.12}
\Pi:V_0^0\to X_{\si}
\end{eqnarray}
 is a linear projection on $X_{\si}$, we have by
(\ref{3.10}) for $\omega \in \mathcal{B}(X_\si)$:
\begin{eqnarray} \label{4.13}
\Pi^{\ast}\nu (\om)=c\int_{B_{\hat{\eps}}} \chi_{\om}(\Pi u)G(du)
=c \int_{\Pi^{-1}\om\cap B_{\hat{\eps}}} G(du)  \nonumber \\
=c \int_{\hat{\Pi}^{-1}(\om\cap \Pi B_{\hat{\eps}})} G(du) =c
\int_{\om\cap \Pi B_{\hat{\eps}}} \Pi^{\ast} G(du).
\end{eqnarray}
We have already mentioned that $\Pi^{\ast} G$ is the Gauss measure
supported on $X_{\si}$ with mathematical expectation $a=0$ and
correlation operator $K_1=\Pi K\Pi^{\ast}$. Note that in
(\ref{4.13}), $\chi_{\om}(v)$ is the characteristic function of
the set $\om$, and
\begin{eqnarray} \label{4.14}
\Pi^{-1}\om=\{x\in V_0^0: \Pi x \in \om \}, \;
\hat{\Pi}^{-1}\om_1=\{x\in V_0^0\cap B_{\hat{\eps}}: \Pi x \in
\om_1 \}.
\end{eqnarray}
Hence to define completely the measure $\Pi^{\ast}\nu(\om)$ from
(\ref{4.13}), we need to calculate $\Pi B_{\hat{\eps}}$. Let us
consider the decomposition (\ref{4.5}) of $V_0^0$. Since
$\Pi:V_0^0\to X_{\si} \subset V_0^0$ is a linear projection, for
each $y \in V_0^0$ decomposed as $y=y'+y''$ with $y'\in X_{\si}$,
$y'' \in X_{\si}^{\perp}$, we obtain
\begin{eqnarray} \label{4.15}
x=\Pi y =y'+A y''\qquad \mbox{where} \quad A=\Pi |_{X^\perp_\si}.
\end{eqnarray}
So
$$
     A: X^\perp_\si \rightarrow X_\si, \qquad A^*: X_\si \rightarrow X_\si^\perp ,
$$
where $A^*$ is the operator adjoint to $A$. (We identify Hilbert
spaces $X_\si, X_\si^\perp$ with their dual spaces.) Note that
$$
B_{\hat{\eps}}=\{ \|y\|^2\leq \hat{\eps}^2\}=\{\|y'\|^2 +\|y''\|^2
\leq \hat{\eps}^2 \}.
$$
Thus by (\ref{4.15}),
\begin{eqnarray} \label{4.16}
\Pi B_{\hat{\eps}} = \{x\in X_{\si}: \text{There exists}\; y''\in
B_{\hat{\eps}}\cap X_{\si}^{\perp} \text{such that}\; \|x-A
y''\|^2+\|y''\|^2 \leq \hat{\eps}^2 \}.
\end{eqnarray}

To make this more precise, we consider the extreme problem
\begin{eqnarray} \label{4.17}
f(y'') \equiv \| x-Ay''\|^2+\|y''\|^2 \to \inf,\;\;\; y''\in
X_{\si}^{\perp}.
\end{eqnarray}
The solution $\hat{y}$ of this problem exists, is unique and
satisfies
$$
(f'(\hat y), h) = 2\{-(x-A\hat y, Ah)+(\hat y,h)\}
  = 2(\hat y+A^{\ast}A\hat y-A^{\ast}x, h)=0, \; \forall h \in X_{\si}^{\perp}.
$$
Since the operator $A^{\ast}A$ is nonnegative, the equation
$$
A^{\ast}A y'' + y''= A^{\ast}x
$$
for each $x\in X_{\si}$ has a unique solution $\hat{y}=\hat{y}(x)
\in X_{\si}^{\perp}$ and it is the solution to the extreme problem
(\ref{4.17}). It is clear that the map $x \mapsto \hat{y}(x):
 X_{\si} \to X_{\si}^{\perp}$ is a bounded linear operator and
 $\hat{y}(x) = (A^{\ast}A+E)^{-1}A^{\ast}x$.
 Thus the definition of $\Pi B_{\hat{\eps}}$ in (\ref{4.16}) can
 be rewritten as follows
\begin{eqnarray} \label{4.18}
\Pi B_{\hat{\eps}} = \{x\in X_{\si}:   \|x-A^{\ast}\hat{y}(x)\|^2
      +\|\hat{y}(x)\|^2 \leq \hat{\eps}^2 \}.
\end{eqnarray}
The  set (\ref{4.18}) is an ellipsoid. Thus (\ref{4.13}) and
(\ref{4.18}) define the measure $\Pi^{\ast}\nu$, and hence the RDS
(\ref{3.1}) is completely defined as well.

 \subsection{Finite-dimensional measure }\label{s4.4}

In order to prove that the RDS (\ref{3.1}) is ergodic, we study a
map of the measure $\Pi^{\ast}\nu$. We introduce the operator of
orthogonal projection $Q$ connected with projection operator $Q_k$
from Subsection \ref{ergodic}
\begin{eqnarray} \label{4.200}
Q: V_0^0 \to X_{\si_k}^{\perp} \equiv X_{\si}^{\perp} \oplus
X_{\si\si_k},
\end{eqnarray}
where subspaces $X_\si^\perp$, $X_{\si\si_k}$ are defined in
(\ref{4.5}),(\ref{4.51}). We have the following property for the
operator $Q$.

\begin{lem} \label{Qlemma}
Let $\Pi, Q$ be projection operators defined in (\ref{2.105}) and
(\ref{4.200}), respectively. Then
\begin{eqnarray} \label{4.21}
Q \Pi = Q\Pi Q.
\end{eqnarray}
\end{lem}

\begin{proof}
Each $x\in V_0^0$ can be decomposed as follows
$$
x=y_1+y_2+y_3
$$
with $y_1\in X_{\si}^{\perp}$, $y_2 \in X_{\si\si_k}$ and $y_3\in
X_{\si_k}$. Then
$$
Q\Pi x = Q(\Pi y_1+y_2+y_3)=Q\Pi y_1+y_2, \quad Q\Pi Q x = Q(\Pi
y_1+y_2)=Q\Pi y_1+y_2.
$$
\end{proof}

The random dynamical system (\ref{3.1}) generates naturally the
measure $\Pi^{\ast}\nu$.  We now study the measure
$Q^{\ast}\Pi^{\ast}\nu$. By (\ref{4.21}) we show that this study
can be reduced to the study of a measure defined on a finite
dimensional space.

\begin{theorem}
Let $Q$ be the orthogonal projection in (\ref{4.200}), $\Pi$ be
the projection in (\ref{2.105}), and $\nu (\om)$ be a probability
measure on $\cB (V_0^0)$. Then
\begin{eqnarray} \label{4.230}
Q^{\ast} \Pi^{\ast} \nu =  Q^{\ast} \Pi^{\ast} Q^{\ast} \nu.
\end{eqnarray}
\end{theorem}

\begin{proof}
By (\ref{3.10}) and (\ref{4.21}), we get
\begin{eqnarray} \label{4.24}
\int f(u)(Q^{\ast} \Pi^{\ast} \nu)(du) =\int f(Q \Pi v)\nu
(dv)          \nonumber \\
 =\int f(Q \Pi Q w) \nu(dw)
 =\int f(u)(Q^{\ast} \Pi^{\ast}  Q^{\ast} \nu)(du)
\end{eqnarray}
This proves the theorem.
\end{proof}

The relation (\ref{4.230}) allows us to reduce our investigation
to the case of measures defined on finite dimensional space. Let
us calculate the measure $Q^{\ast} \nu$. Taking into account the
definitions (\ref{3.11}) of $\nu (du)$ and (\ref{4.200}) of $Q$,
we get, analogous to (\ref{4.13}), that for each $\Gamma \in
\cB(X_{\si_k}^{\perp})$,
\begin{eqnarray} \label{4.250}
Q^{\ast}\nu (\Gamma)=c \int_{\Gamma \cap QB_{\hat{\eps}}}
Q^{\ast}G(du),
\end{eqnarray}
where $Q^{\ast}G$ is the Gauss measure supported on
$X_{\si_k}^{\perp}$, with   mathematical expectation zero and
correlation operator $QKQ$. We make the following identification
taking into account identity in (\ref{4.200}): using in
$X^\perp_{\si_k}$ an orthonormal basis $\{ e_j,\, j=1, \dots
,n\}\, (n=n_k)$ introduced in Subsection \ref{s4.1}, we see can
write
\begin{equation}\label{4.231}
y=\sum_{j=1}^ny_je_j\in X^\perp_{\si_k},\quad
u=\sum_{j=1}^mu_je_j\in X^\perp_{\si},\quad
v=\sum_{j=m+1}^nv_je_j\in X_{\si \si_k},
\end{equation}
and take the following identifications:
\begin{equation}\label{4.232}
X^\perp_{\si_k}\cong \bR^n=\{ \vec{y}=(y_1,\dots ,y_n)\},\;
X^\perp_{\si}\cong \bR^m=\{ \vec{u}\}, \; X_{\si \si_k}\cong
\bR^{n-m}=\{ \vec{v}\}.
\end{equation}
 We restrict correlation operator $QKQ$ on
$X_{\si_k}^{\perp}$.  Then $QKQ$ can be regarded as
   a $n \times n$ matrix. By (\ref{3.7}) this matrix is non-degenerate because for each
$0\ne u\in X^\perp_{\si_k}\; (u,QKQu)=(Qu,KQu)=(u,Ku)>0$ and
therefore $\ker QKQ=0$ \footnote{Actually we can weaken the first
condition in (\ref{3.7}) assuming that $K^*=K\ge 0$ and $\ker
QKQ=0$ where $Q$ is orthogonal projector on $X^\perp_{\si_k}$ with
big enough $\si_k$} We denote $ \hat{K} = (QKQ)^{-1}$. By
(\ref{3.8})-(\ref{3.9}), we conclude that
\begin{eqnarray} \label{4.26}
Q^{\ast}G(dy) := \hat{G}(dy)=g(y)dy,
\end{eqnarray}
where $g(y)= \frac{\det \hat{K}}{(2\pi)^{\frac{n}2}}
e^{-\frac12(\hat{K}y,y)}$. Note that
\begin{eqnarray} \label{4.26'}
QB_{\hat{\eps}}=\{y=\sum_{j=1}^ny_je_j \in X^\perp_{\si_k}:
\sum_{j=1}^n y_j^2 \leq \eps^2 \} := {B}.
\end{eqnarray}
We hence obtain
\begin{eqnarray} \label{4.27}
\hat{\nu}(dy) := Q^{\ast}\nu(dy)=\hat{c}\hat{\chi}_{\hat{\eps}}(y)
\hat{G}(dy),
\end{eqnarray}
where $\hat{c} =(\int_{B} \hat{G}(dx))^{-1}$ and
$\hat{\chi}_{\hat{\eps}}$ is the characteristic function of the
ball ${B}$.

We  introduce the projection operator
\begin{eqnarray}\label{5.2}
  \pi=Q\Pi: X_{\si_k}^{\perp} \to X_{\si\si_k},
\end{eqnarray}
In virtue of (\ref{4.27}), for each $\om \in \cB
(X^\perp_{\si_k})$,
\begin{eqnarray} \label{4.29}
Q^{\ast}\Pi^{\ast}Q^{\ast}\nu(\om)= \pi^{\ast} \hat{\nu}(\om).
\end{eqnarray}


\section{Density $P(x)$ of the measure $\pi^{\ast} \hat{\nu}$ }

\subsection{Preliminaries}

Recall that the space $X_{\si_k}^{\perp}$ admits the orthogonal
decomposition
\begin{eqnarray}\label{5.1}
  X_{\si_k}^{\perp}= X_{\si}^{\perp} \oplus  X_{\si\si_k}.
\end{eqnarray}
Below in order to emphasize belonging $u\in X^\perp_\si , v\in
X_{\si \si_k}$ we write $u\oplus v$ instead of $u+v$.

The projection operator $\pi$ defined in (\ref{5.2})
 can be represented as follows
\begin{eqnarray}\label{5.3}
  \pi = (\alpha, E),
\end{eqnarray}
where $E$ is the identity operator in $X_{\si\si_k}$ and $\alpha:
X_{\si}^{\perp} \to X_{\si\si_k}$. For  $x\in X_{\si\si_k}$ we
denote by $\pi_x$,
 the affine plane in $X_{\si_k}^{\perp }$:
\begin{eqnarray}\label{5.4}
  \pi_x = \pi^{-1} x
  =\{y\equiv u\oplus v \in X_{\si}^{\perp} \oplus
  X_{\si\si_k}=X_{\si_k}^{\perp}: \alpha u+v=x \}.
\end{eqnarray}
In particular, when $x=0$,
\begin{eqnarray}\label{5.5}
  \pi_0 = \pi^{-1} 0
  =\{y\equiv u\oplus v \in X_{\si}^{\perp} \oplus
  X_{\si\si_k}=X_{\si_k}^{\perp}: \alpha u+v=0 \}.
\end{eqnarray}
Since $B$ is the support of the measure $\hat{\nu}(dy)\equiv
\hat{\nu}(du, dv)$ defined in (\ref{4.27}), the ellipsoid $\pi B$
is the support of the measure $\pi^{\ast}\hat{\nu}$.

For each $f\in C(\pi B)$ we have
\begin{eqnarray}\label{5.6}
  \int_{\pi B}f(x)\pi^{\ast}\hat{\nu}(dx)=\int_B f(\alpha
  u+v)\hat{\nu}(du,dv)=\int_{\pi B}f(x)dx \int_{\pi_x\cap
  B}\Gamma (w, x)d w,
\end{eqnarray}
where the first equality is via the definition of the measure
$\pi^{\ast}\hat{\nu}$ and the second equality follows from the
change of variables
\begin{eqnarray}\label{5.7}
  w=u-\alpha u\in \pi_0,\; x=\alpha u+v\in X_{\si_k}^{\perp}.
\end{eqnarray}
The calculation of $\Gamma(w,x)$ will be done later. The formulae
(\ref{5.6}) gives the expression of the density $P(x)$ for the
measure $\pi^{\ast}\hat{\nu}(dx)$:
\begin{eqnarray}\label{5.8}
  P(x)dx = \pi^{\ast}\hat{\nu}(dx), \; \text{where}\; P(x)=\int_{\pi_x\cap
  B}\Gamma (w, x)d w.
\end{eqnarray}


\subsection{Change of variables $(w,x) \to (u,v)$}

In order to calculate the kernel functional $\Gamma (w, x)$ in
(\ref{5.8}), we need to consider the following change of
variables, i.e., the inverse of (\ref{5.7})
\begin{eqnarray}\label{5.9}
  u=u(w,x), \; v=v(w,x).
\end{eqnarray}
We introduce an orthonormal basis $\{b_j, j=1, \cdots, n\}$ in
$X_{\si_k}^{\perp}$. Let
\begin{eqnarray}\label{5.10}
\{b_j, j=1, \cdots, m \} \; \text{be the orthonormal basis of}\;
X_{\si}^{\perp}  \subset  X_{\si_k}^{\perp},
\end{eqnarray}
composed of eigenvectors of the operator $E+\a^{\ast}\a:
X_{\si}^{\perp} \to X_{\si}^{\perp}$. Suppose that $1\leq \mu_1,
\cdots, 1\leq \mu_m$ are eigenvalues corresponding to eigenvectors
$b_1, \cdots, b_m$. Assume that
\begin{eqnarray}\label{5.11}
 \mu_1>1, \cdots, \mu_s>1, \mu_{s+1}=\cdots =\mu_m =1.
\end{eqnarray}

\begin{lem} \label{lemma5.1}
The following statements hold:

 (i) $b_j \in \ker \a, \; j=s+1,\cdots, m$;

(ii) $\{ \a b_j, \; j=1, \cdots, s\}$ form an orthogonal basis in
$Im \a \subset X_{\si\si_k}$.
\end{lem}

\begin{proof}
(i) For $j=s+1,\cdots, m$, $(E+\a^{\ast}\a)b_j = b_j $ iff
$\a^{\ast}\a b_j =0$. Since $\ker \a^{\ast} \perp Im\a$, $\a b_j
\neq 0$ implies $\a^{\ast}\a b_j \neq 0$. Hence $b_j \in \ker \a$.

(ii) Note that for $i, j=1, \cdots, s$, $(\a b_i, \a
b_j)=(\a^{\ast}\a b_i, b_j)=(\mu_i-1)(b_i, b_j)=0$ and $\mu_i \neq
1$. Hence $\a b_i \perp \a b_j$. If $\gamma \in Im\a$ then for a
certain $b\in X_{\si}^{\perp}$, we have $\gamma = \a b= \a
\sum_{j=1}^m c_jb_j =\sum_{j=1}^sc_j\a b_j$. Therefore $\{\a b_j,
j=1, \cdots, s\}$ form a basis in $Im \a$.
\end{proof}

Since $\ker \a^{\ast} \oplus Im \a =X_{\si\si_k}$ and     $\dim
Im\a = s$ due to Lemma \ref{5.1}, we see that $\dim \ker \a^{\ast}
= n-m-s$. Let $\{b_{m+s+1}, \cdots, b_n\}$ be an orthonormal basis
for $\ker \a^{\ast}$. Then by Lemma \ref{5.1} again, the vectors
\begin{eqnarray}\label{5.12}
b_{j+m}=\frac{\a b_j}{\|\a b_j\|}, j=1, \cdots, s,\;\;  b_{m+s+1},
\cdots, b_n
\end{eqnarray}
form an orthonormal basis in $X_{\si\si_k}$. In virtue of
(\ref{5.10})-(\ref{5.12}),
\begin{eqnarray}\label{5.13}
\text{Vectors} \; b_1, \cdots, b_n \; \text{form an orthonormal
basis of}\; X_{\si_k}^{\perp}.
\end{eqnarray}
On the plane $\pi_0$ defined in (\ref{5.5}), we consider the
vectors
\begin{eqnarray}\label{5.15}
\theta_j=\frac{b_j\oplus (-\a b_j)}{(1+ \|\a b_j\|^2)^{\frac12}},
\; j=1,\cdots,m.
\end{eqnarray}

\begin{lem} \label{lemma5.2}
Vectors (\ref{5.15}) form an orthonormal basis of the plane
$\pi_0$.
\end{lem}

\begin{proof}
By the definitions (\ref{5.3}), (\ref{5.5}) and (\ref{5.15}), we
see that $\theta_j \in \pi_0$, as
$$
\pi \theta_j =(1+ \|\a b_j\|^2)^{-\frac12}[\a b_j-\a b_j]=0.
$$
By Lemma \ref{5.1}, for $i \neq j$,
\begin{align}   \label{5.15'}
(\theta_i, \theta_j)_{\pi_0}
 &= (\frac{b_i\oplus (-\a
b_i)}{\sqrt{1+ \|\a b_i\|^2}},  \frac{b_j\oplus (-\a
b_j)}{\sqrt{1+
\|\a b_j\|^2}} )_{X^\perp_{\si_k}}   \nonumber \\
 & =  (1+ \|\a
b_i\|^2)^{-\frac12} (1+ \|\a b_j\|^2)^{-\frac12}
[(b_i,b_j)_{X_{\si}^{\perp}} + (-\a b_i, -\a
b_j)_{X_{\si\si_k}}]   \nonumber\\
 &=  0
\end{align}
and therefore the system (\ref{5.15}) is orthonormal. As the rank
of the matrix for (\ref{5.3}) equals to $n-m$,
$$
\dim \pi_0=\dim \ker \pi \equiv \dim \ker (\a, E)=m.
$$
Hence the system (\ref{5.15}) forms an orthonormal basis of the
plane $\pi_0$.
\end{proof}

Define
\begin{eqnarray}\label{5.16}
 \theta_j = b_j, \; j=m+1, \cdots, n.
\end{eqnarray}
Then the vectors
\begin{eqnarray}\label{5.17}
 \theta_j, \; j=1, \cdots, n.
\end{eqnarray}
defined in (\ref{5.15}) and (\ref{5.16}) form a basis of
$X_{\si_k}^{\perp}$.

Let $R=(R_{ij})$ be the $n \times n$ matrix with components
$R_{ij}$ defined as follows.
\begin{align}\label{5.18}
 R_{ii} & = (1+ \|\a b_i\|^2)^{-\frac12};
 & R_{i,m+1}=\frac{-\alpha }{\sqrt{1+ \|\a b_i\|^2}},\; i=1,\cdots, s; \nonumber \\
R_{ii} & =1, \; i=s+1,\cdots, n;
 & R_{ij}=0   \; \text{for other $i,j$}.
\end{align}
It can be checked directly that the matrix $R$ transforms the
basis $\vec{b}=(b_1, \cdots, b_n)$ to the basis
$\vec{\theta}=(\theta_1, \cdots, \theta_n)$:
\begin{eqnarray}\label{5.19}
 \vec{\theta} =R \vec{b}.
\end{eqnarray}

Now we can calculate the change of variables (\ref{5.9}), the
inverse of (\ref{5.7}). Let $y\in X_{\si_k}^{\perp}$ admits
decompositions $y=u+v=w+x$ with $u\in X_{\si}^{\perp}$, $v\in
X_{\si\si_k}$, $w\in \pi_0$ and $x\in X_{\si\si_k}$. Define
\begin{eqnarray}\label{5.20}
 y=\sum_{j=1}^n y_jb_j=\sum_{j=1}^n z_j\theta_j, \nonumber \\
 u=\sum_{j=1}^mu_jb_j,       \; v=\sum_{j=m+1}^n v_{j-m}b_j,     \nonumber \\
 w=\sum_{j=1}^m w_j\theta_j, \; x=\sum_{j=m+1}^nx_{j-m}\theta_j.
\end{eqnarray}

We introduce notations
\begin{eqnarray}\label{5.21}
  \vec{y}=(y_1,\cdots,y_n),
  \vec{z}=(z_1,\cdots,z_n),                       \nonumber\\
  \vec{u}=(u_j\equiv y_j,j=1,\cdots,m),
  \vec{v}=(v_j\equiv y_{j+m}, j=1,\cdots, n-m),  \nonumber\\
  \vec{w}=(w_j\equiv z_j, j=1,\cdots,m),
  \vec{x}=(x_j\equiv z_{j+m}, j=1,\cdots,n-m).
\end{eqnarray}
Then the change of variables (\ref{5.9}) is rewritten as
\begin{eqnarray}\label{5.22}
 \vec{u}=u(\vec{w},\vec{x}), \; \vec{v}=v(\vec{w}, \vec{x}),
\end{eqnarray}
or in more compact form
\begin{eqnarray}\label{5.23}
 \vec{y}=y(\vec{z}).
\end{eqnarray}

\begin{theorem} \label{transform}
The transformation (\ref{5.23}) can be calculated as follows
\begin{eqnarray}\label{5.24}
 \vec{y}=y(\vec{z})=R^{\ast} \vec{z},
\end{eqnarray}
where $R^{\ast}$ is the conjugate matrix of the matrix $R$ defined
in (\ref{5.18})-(\ref{5.19}). Note that the matrix $R$ transforms
the basis $\{b_j\}$ to the basis $\{\theta_j\}$.
\end{theorem}

\begin{proof}
In fact, the relation (\ref{5.24}) follows from
(\ref{5.19})-(\ref{5.20}).
\end{proof}

Using (\ref{5.24}) and (\ref{5.18}), we obtain the Jacobian of the
transformation (\ref{5.24}):
\begin{eqnarray}\label{5.25}
 J= \det D\vec{y}/D\vec{z} =\det R^{\ast}
 =\prod_{i=1}^s (1+ \|\a b_i\|^2)^{-\frac12}.
\end{eqnarray}
We see that the Jacobian $J$ depends only on the operator $\a$.
Now we make more precise the expression for density $P(x)$ in
(\ref{5.8}). Note that in (\ref{5.6}) we made just the change of
variables (\ref{5.23}) or the equivalent (\ref{5.24}). Taking into
account of the definition $\vec{y}=(\vec{u}, \vec{v})$ and
$\vec{z}= (\vec{w}, \vec{x})$, and using the facts $(\ref{4.26}),
(\ref{4.27}), (\ref{5.24})$ and $(\ref{5.25})$, we obtain that the
integrand $\Gamma(w,x)$ in the expression of density $P(x)$ in
(\ref{5.8}):
\begin{eqnarray}\label{5.26}
\Gamma(w,x) = \prod_{i=1}^s (1+ \|\a b_i\|^2)^{\frac12} \det
\hat{K} (2\pi)^{-n/2} \exp[-\frac12 (\hat{K}QR^{\ast}\vec{z},
QR^{\ast}\vec{z})],
\end{eqnarray}
with $\vec{z}= (w, x)$. We suppress the arrow $\to$ on top of $w$
and $x$ here.


\section{Smoothness of the density $P(x)$}

Formulas (\ref{5.6}) and (\ref{5.8}) imply that the density $P(x)$
is supported in the ellipsoid  $\pi B\subset X_{\si\si_k}$. So
$P(x) \equiv 0$ for $x\in X_{\si\si_k}\setminus \pi B$. We now
investigate the smoothness of $P(x)$ for $x\in \pd (\pi B)$ and
for $x\in \text{Int} (\pi B)$, respectively.


\subsection{Smoothness of the density $P(x)$ on boundary $\pd (\pi B)$}

It is clear that the set $B\cap \pi_x$ is $\emptyset$ if $x \notin
\pi B$, it is a single point if $x \in \pd (\pi B)$, and it is a
ball in the $n-m$ dimensional plane $\pi_x$ if $x \in \text{Int}
(\pi B)$. We first calculate the center and radius of this ball.

\begin{lem} \label{lemma6.1}
Let the ball $B$ be defined in (\ref{4.26'}),
 $x \in \text{Int} (\pi B)$, and the plane
$\pi_x$ be defined in (\ref{5.4}). Then the center of the ball
$B\cap \pi_x$ is
\begin{eqnarray}\label{6.1}
\tilde{w}=\hat{u}\oplus \hat{v}=
[(E+\a^{\ast}\a)^{-1}\a^{\ast}x]\oplus
[x-\a(E+\a^{\ast}\a)^{-1}\a^{\ast}x ]
\end{eqnarray}
and the radius of the ball $B\cap \pi_x$ is
\begin{eqnarray}\label{6.2}
r= (\eps^2-
\|(E+\a^{\ast}\a)^{-1}\a^{\ast}x\|^2_{X_{\si}^{\perp}}-
\|x-\a(E+\a^{\ast}\a)^{-1}\a^{\ast}x\|^2_{X_{\si\si_k}}
)^{\frac12},
\end{eqnarray}
where, recall that, $\eps$ is the radius of the ball $B$.
\end{lem}

\begin{proof}
Evidently, the center $\{\hat{u},\hat{v} \}$ is the solution of
the extreme problem
\begin{eqnarray}\label{6.3}
 \|u\|^2_{X_{\si}^{\perp}}+\|v\|^2_{X_{\si\si_k}} \to \inf,\;
 \{u,v\} \in \pi_x.
\end{eqnarray}
By definition (\ref{5.4}), $\{u,v\} \in \pi_x$ iff $\a u+v=x$,
i.e., $v=x-\a u$. Substituting this into (\ref{6.3}) and solving
the extreme problem, we obtain the solution (\ref{6.1}). The
radius $r$ follows from the Pythagoras theorem.
\end{proof}

Let us consider the following extreme problem: Given $x\in \pi B$,
find $h\in X_{\si\si_k}$ such that
\begin{align}
 F(h)    &   =\|(E+\a^{\ast}\a)^{-1}\a^{\ast}(x+h)\|^2
 +\|(E-\a(E+\a^{\ast}\a)^{-1}\a^{\ast})(x+h)\|^2 \to \inf, \label{6.4}\\
 \|h\|^2 & =\gamma_0^2,            \label{6.5}
\end{align}
with $\gamma_0>0$ a given sufficiently small parameter. Recall
that each $x\in X_{\si\si_k}$ admits the decomposition
\begin{align} \label{6.6}
x=x_0+\sum_{j=1}^s x_j\a b_j, \; x_0\in \ker \a^{\ast},\; x_j\in
\bR,
\end{align}
where $\{b_j\}$ is the basis in (\ref{5.10}). Note also that, by
Lemma \ref{lemma5.1}(ii), $\{\a b_j, j=1, \cdots, s\}$ is a basis
of $Im \a$.

\begin{lem} \label{lemma6.2}
Suppose that $x\in \pi B$ and it has decomposition (\ref{6.6}). If
$\gamma_0>0$ is small enough, then there exists a unique solution
$\hat{h}$ of the extreme problem (\ref{6.4})-(\ref{6.5}). The
solution $\hat{h}$ is determined by
\begin{align} \label{6.7}
\hat{h}= h_0+\sum_{j=1}^s h_j\a b_j, \; h_0\in \ker \a^{\ast},\;
h_j\in \bR,
\end{align}
where
\begin{align} \label{6.8}
h_0=-\frac{x_0}{1+\l (\gamma_0)}, \; h_j=-\frac{x_j}{1+\l
(\gamma_0)\mu_j}, \; j=1,\cdots, s
\end{align}
and $x_0, x_j, j=1,\dots ,s$ are defined in (\ref{6.6}), $\mu_j>1,
\; j=1,\cdots, s$ are eigenvalues (\ref{5.11}) of the operator
$E+\a^{\ast}\a$, and $\l(\gamma_0)$ is the unique solution of the
equation
\begin{align} \label{6.9}
 \frac{\|x_0\|^2}{(1+\l)^2}
 +\sum_{j=1}^s \frac{x_j^2 \|\a
 b_j\|^2}{(1+\l\mu_j)^2}=\gamma_0^2.
\end{align}
\end{lem}

\begin{proof}
The existence of a solution of the finite-dimensional problem
(\ref{6.4})-(\ref{6.5}) is evident. We now prove the uniqueness of
this solution $h$.

Let $\cL (h,\l)=F(h)+\l (\|h\|^2-\gamma_0^2)$ be the Lagrange
function for the extreme problem (\ref{6.4})-(\ref{6.5}). By the
Lagrange principle, if $h$ is a solution of this problem, then
there exists $\l \in \bR$ such that
\begin{align} \label{6.10}
 (\cL'_h(h, \l), \delta)=(F'(h), \delta)+ 2\l (h,\delta)=0,\;
 \forall \delta \in  X_{\si\si_k}.
\end{align}
Substitution of expression $F(h)$ in (\ref{6.4}) into (\ref{6.10})
yields
\begin{align*}
((E+\a^{\ast}\a)^{-1}\a^{\ast}(x+h),
(E+\a^{\ast}\a)^{-1}\a^{\ast}\delta)   \nonumber \\
 +((E-\a(E+\a^{\ast}\a)^{-1}\a^{\ast})(x+h),
 (E-\a(E+\a^{\ast}\a)^{-1}\a^{\ast})\delta)
 +\l (h, \delta) = 0.
\end{align*}
We transform operators from right multiplies in scalar products to
the left multipliers. Noting that
\begin{align} \label{6.10'}
  \a(E+\a^{\ast}\a)^{-2}\a^{\ast}+(E-\a(E+\a^{\ast}\a)^{-1}\a^{\ast})^2
  =E-\a(E+\a^{\ast}\a)^{-1}\a^{\ast},
\end{align}
we obtain equations for $h$ and $\l$:
\begin{align} \label{6.11}
(E-\a(E+\a^{\ast}\a)^{-1}\a^{\ast}) (x+h) + \l h = 0.
\end{align}
Substituting the decompositions (\ref{6.6})-(\ref{6.7}) for $x$
and $h$ into (\ref{6.11}), we arrive at the following system of
equations
\begin{equation}\label{6.12}
x_0 + h_0(1+\l) = 0,
\end{equation}
\begin{equation}\label{6.13}
[1-\frac{\mu_j-1}{\mu_j}](x_j+h_j)+\l h_j = 0,
\end{equation}
and therefore
\begin{equation} \label{6.14}
h_0=-\frac{x_0}{1+\l}, \; h_j=-\frac{x_j}{1+\l\mu_j}.
\end{equation}
By (\ref{6.5}),(\ref{6.7}) and (\ref{6.14}), we finally get
\begin{equation} \label{6.15}
 r(\l) \equiv \|h\|^2=\|h_0\|^2+\sum_{j=1}^s h_j^2\|\a b_j\|^2
 =   \frac{\|x_0\|^2}{(1+\l)^2}
 +\sum_{j=1}^s \frac{x_j^2 \|\a
 b_j\|^2}{(1+\l\mu_j)^2}=\gamma_0^2.
\end{equation}
Since $r(0)=\|x\|^2$, $r(\l)\to 0$ as $\l \to \infty$ and
$r'(\l)<0$, there exists a unique solution $\l =\l(\gamma_0)$ if
$\gamma_0^2 < \|x\|^2$.
\end{proof}

\begin{theorem} \label{theorem6.1}
Let $x \in \partial \pi B$, $x+h\in \text{Int}(\pi B)$ and $\|h\|
< \|x\|$. Then $P(x)=0$ and for some positive constant $c$,
\begin{equation} \label{6.16}
  |P(x+h)|\leq c \|h\|^{\frac{m}2}\;\;
 \text{as}\;\;  \|h\| \to 0.
\end{equation}
\end{theorem}

\begin{proof}
If $x\in \pd(\pi B)$, then
$$
B\cap \pi_x = \hat{u}\oplus \hat{v},
$$
where the point $\hat{u}\oplus \hat{v} $ is defined  in
(\ref{6.1}). Hence by  (\ref{5.8}), $P(x)=0$. Denote
$\|h\|=\gamma_0$ and take the solution $\hat{h}$ of the problem
(\ref{6.4})-(\ref{6.5}). By (\ref{6.2}),(\ref{6.4}) and
(\ref{6.5}), the radius $r(x+\hat{h})$ of the ball $B \cap
\pi_{x+\hat{h}}$ is maximal in the sets of radius $r(x+h)$ of $B
\cap \pi_{x+h}$ corresponding to vectors $h$ such that
$\|h\|=\gamma_0$. We calculate   $r(x+\hat{h})$.

Since $x\in \pd(\pi B)$, we have
$\|\hat{u}(x)\|^2+\|\hat{v}(x)\|^2=\eps^2$ with $\{\hat{u},
\hat{v} \}$   defined  in (\ref{6.1}). In (\ref{6.2}),   taking
$x\to x+\hat{h}$ and substituting
$\eps^2=\|\hat{u}(x)\|^2+\|\hat{v}(x)\|^2$, and taking into
account of (\ref{6.6})-(\ref{6.8}), we get
\begin{align} \label{6.17}
r^2(x+\hat{h}) & =-( (E+\a^{\ast}\a)^{-1}\a^{\ast}\hat{h},
(E+\a^{\ast}\a)^{-1}\a^{\ast}(2x+\hat{h}) )   \nonumber \\
&-((E-\a(E+\a^{\ast}\a)^{-1}\a^{\ast})\hat{h},
(E-\a(E+\a^{\ast}\a)^{-1}\a^{\ast})(2x+\hat{h}) ) \nonumber \\
&=( (E+\a^{\ast}\a)^{-1}\a^{\ast} \sum_{j=1}^s \frac{x_j\a
 b_j}{1+\l\mu_j},
(E+\a^{\ast}\a)^{-1}\a^{\ast}\sum_{l=1}^sx_l\a b_l(2-\frac1{1+\l\mu_j})  )   \nonumber \\
&+
(\frac{x_0}{1+\l}+\sum_{j=1}^s(E-\a(E+\a^{\ast}\a)^{-1}\a^{\ast})\a
b_j\frac{x_j}{1+\l\mu_j},             \nonumber \\
& x_0(2-\frac1{1+\l})+\sum_{l=1}^s
(E-\a(E+\a^{\ast}\a)^{-1}\a^{\ast})\a b_lx_l(2-\frac1{1+\l\mu_j}) )    \nonumber \\
& =\sum_{j,l=1}^s \frac{x_j(\mu_j-1)}{\mu_j(1+\l \mu_j)}
\frac{x_l(\mu_l-1)}{\mu_l}(2-\frac1{1+\l\mu_l}) (b_j,b_l)
+\frac{\|x_0\|^2(1+2\l)}{(1+\l)^2}    \nonumber \\
& +\sum_{j,l=1}^s x_jx_l(\a b_j, \a
b_l)(1-\frac{\mu_j-1}{\mu_j})\frac1{1+\l\mu_j}(1-\frac{\mu_l-1}{\mu_l})(2-\frac1{1+\l\mu_l})
\nonumber \\
&= \sum_{j=1}^s\frac{x_j^2(\mu_j-1)^2(1+2\l \mu_j)}{\mu_j^2(1+\l
\mu_j)^2} + \frac{\|x_0\|^2(1+2\l)}{(1+\l)^2} + \sum_{j=1}^s x_j^2
\|\a b_j\|^2 \frac{1+2\l\mu_j}{\mu_j^2(1+\l \mu_j)^2}.
\end{align}
We rewrite this as follows:
\begin{align} \label{6.18}
r^2(x+\hat{h}) =\frac{\|x_0\|^2(1+2\l)}{(1+\l)^2}
+\sum_{j=1}^s\frac{x_j^2\|\a b_j\|^2 }{(1+\l \mu_j)^2} \frac{1+2\l
\mu_j}{\mu_j^2} (1+\frac{(\mu_j-1)^2}{\|\a b_j\|^2}).
\end{align}

There exist constants $0<c_1<c_2/2$ such that for each $j=1,
\cdots, s$ and for every $\l >0$, we have
\begin{align} \label{6.19}
  c_1(1+\l) \leq \frac{1+2\l\mu_j}{\mu_j^2} (1+\frac{(\mu_j-1)^2}{\|\a b_j\|^2})
  \leq \frac{c_2}{2} (1+\l).
\end{align}
Comparing (\ref{6.9}) and (\ref{6.18}), we thus get
\begin{align} \label{6.20}
  c_1\gamma_0^2(1+\l) \leq r^2(x+\hat{h})  \leq c_2 \gamma_0^2 (1+\l).
\end{align}

Let $A_1(x)=(\|x_0\|^2+\sum_{j=1}^s\frac{x_j^2\|\a b_j\|^2
}{\mu_j^2})^{\frac12}$. Then using (\ref{6.9}) we get
\begin{align*}
   \gamma_0^2 \geq  \frac{\|x_0\|^2}{(1+\l)^2}
+\sum_{j=1}^s\frac{x_j^2\|\a b_j\|^2 }{\mu_j^2(1+\l)^2} =
\frac{A_1(x)^2}{(1+\l(\gamma_0))^2}.
\end{align*}
and therefore
\begin{align} \label{6.21}
   1+\l(\gamma_0) \geq \frac{A_1(x)}{\gamma_0}.
\end{align}

Let now $A_2(x)=(\|x_0\|^2+\sum_{j=1}^sx_j^2\|\a
b_j\|^2)^{\frac12}$. Then by (\ref{6.9}) we obtain
$$
\gamma_0^2 \leq \frac{A_2(x)^2}{(1+\l)^2}.
$$
Hence,
\begin{align} \label{6.22}
   1+\l(\gamma_0) \leq \frac{A_2(x)}{\gamma_0}.
\end{align}
Substituting (\ref{6.21}), (\ref{6.22}) into (\ref{6.20}), we get
\begin{align} \label{6.23}
  c_1 A_1(x)\gamma_0 \leq r^2(x+\hat{h})\leq c_2 A_2(x)\gamma_0.
\end{align}
By (\ref{5.26}), there exists constants $0<\hat{c}_1 <\hat{c}_2$
such that
\begin{align} \label{6.24}
  \hat{c}_1 \leq \Gamma (w,x) \leq   \hat{c}_2,
\end{align}
for each $(w,x)=\vec{z}$ such that $\|QR^{\ast}\vec{z}\|^2 \equiv
\|\vec{y}\|^2 \leq \eps^2$, i.e., on the ball $B$. Hence, by
(\ref{5.8}),
\begin{align} \label{6.25}
  \hat{c}_1 V(\pi_{x+h}\cap B) \leq P(x+h)=\int_{\pi_{x+h}\cap B}\Gamma(w,x)dw
    \leq   \hat{c}_2 V(\pi_{x+h}\cap B) ,
\end{align}
where $V(\pi_{x+h}\cap B)$ is the volume of the ball
$\pi_{x+h}\cap B$.

Note that $r(x+h) \leq r(x+\hat{h})$ for each $h\in X_{\si\si_k}$
such that $\|h\|=\|\hat{h}\|=\gamma_0$ and $x+h\in \text{Int} (\pi
B)$. Thus (\ref{6.23}),(\ref{6.25}) and the fact that
$V(\pi_{x+h}\cap B) = c_m r(x+h)^m$ (where $c_m$ is the volume of
the unit ball in $\bR^m$) implies (\ref{6.16}).
\end{proof}

{\bf Remark.} The inequalities in (\ref{6.25}) imply that for each
$x\in \pd (\pi B)$ and for $\hat{h}$ defined in the proof of
Theorem \ref{theorem6.1}, the following estimate holds
$$
P(x+\hat{h}) \geq c \|\hat{h}\|^{\frac{m}2}.
$$
This inequality means that if $m=1$, $P(x)$ is not differentiable
at $x$ on the boundary $\pd (\pi B)$.


\subsection{Smoothness of the density $P(x)$ in the interior $\text{Int} (\pi B)$}


We first reformulate the definition of $P(x)$ in (\ref{5.8}). By
(\ref{5.4}) and (\ref{5.5}), for each $x\in X_{\si\si_k}$, we see
that
$$
\pi_x = 0\oplus x + \pi_0, \;\text{with}\; \pi_0=\{u\oplus (-\a
u), u\in X_{\si}^{\perp}\} \subset X_{\si}^{\perp} \oplus
X_{\si\si_k}.
$$
Therefore the ball
$$
\pi_x \cap B =\{w=u\oplus v \in \pi_x: \|w-\hat{w}\|^2 \leq r^2\},
$$
where the center $\hat{w}=\hat{w}(x)$ and the radius $r=r(x)$ are
defined in (\ref{6.1})-(\ref{6.2}), can be rewritten as
$$
\pi_x \cap B =0\oplus x+ \cD(x).
$$
Here
\begin{align} \label{6.27}
\cD(x) = \{w=u\oplus (-\a u)\in\pi_0: \|w-\hat{w}(x)\|^2 \leq
r^2(x) \}
\end{align}
with $r(x)$ being defined in (\ref{6.2}) and
\begin{align} \label{6.28}
\hat{w}(x)=[(E+\a^{\ast}\a)^{-1}\a^{\ast}x]\oplus
[-\a(E+\a^{\ast}\a)^{-1}\a^{\ast}x].
\end{align}
We now decompose $w$ and $\hat{w}$ in the basis $\theta_j$ of
$\pi_0$, defined in (\ref{5.15}):
$$
w=\sum_{j=1}^m w_j\theta_j, \; \hat{w}(x)=\sum_{j=1}^m
\hat{w}_j\theta_j,
$$
and consider the integral
\begin{align} \label{6.29}
P(x)=\int_{B\cap\pi_x} \Gamma(w,x)dw =\int_{\cD(x)} \Gamma(w,x)dw,
\end{align}
with $w=(w_1,\cdots, w_m) \in \bR^m$.

Let $x\in \text{Int}(\pi B)$. To investigate the smoothness of
$P(x)$ at $x$, we consider the difference
\begin{align} \label{6.30}
P(x+h)-P(x)=\int_{\cD(x+h)}\Gamma(w,x+h)dw
-\int_{\cD(x)}\Gamma(w,x)dw.
\end{align}
By (\ref{5.26}), the function $\Gamma(w,x)$ is infinitely
differentiable in $x$ and in $w$. By Taylor expansion
\begin{align} \label{6.31}
 \Gamma(w,x+h)= \Gamma(w,x) +(\Gamma_x'(w,x), h)+ o(w,x,h),
\end{align}
where
\begin{align} \label{6.32}
   |o(w,x,h)| \leq c \|h\|^2
\end{align}
with $c$ depending on $(w,x)$. As usual, in this case, we denote
$o(w,x,h)$ as $O(\|h\|^2)$. Inserting (\ref{6.31}) into
(\ref{6.30}), we conclude that
\begin{align} \label{6.33}
P(x+h)-P(x) &=\int_{\cD(x+h)\setminus \cD(x)}\Gamma(w,x)dw
-\int_{\cD(x)\setminus \cD(x+h)}\Gamma(w,x)dw  \nonumber  \\
& + \int_{\cD(x+h)}(\Gamma_x'(w,x), h)dw + o(\|h\|^2).
\end{align}

We now calculate the G\^{a}teaux  derivative of $P(x)$ at $x$.
Denote $e=h/\|h\|$ and $ \l =\|h\|$. Divide (\ref{6.33}) by $\l$
and take limit as $\l \to 0$. For the third term in the right hand
side, we have
\begin{align} \label{6.34}
\lim_{\l\to 0} \frac1{\l} \int_{\cD(x+\l e)}(\Gamma_x'(w,x), \l
e)dw = \int_{\cD(x)}(\Gamma_x'(w,x), e)dw.
\end{align}

To find the similar limit for the first and second terms in the
right hand side, we introduce a kind of ``polar" coordinates in
the sets $\cD(x+h)\setminus \cD(x)$ and $ \cD(x)\setminus
\cD(x+h)$ when they are not empty. Suppose that
\begin{align} \label{6.35}
\cD(x+h)\setminus \cD(x) \neq \emptyset.
\end{align}

Let $b$ be the running point on the part of the sphere $\pd
B(x+h)$ which is the part of boundary for the set (\ref{6.35}).
Define
\begin{align} \label{6.36}
a=\overline{w(x)b} \cap \pd B(x),
\end{align}
where $\overline{w(x)b}$ is the vector with end points $w(x)$ and
$b$. Denote by  $\psi$  the magnitude of the angle $\angle b w(x)
w(x+h)$. By the cosine theorem in the triangle $\triangle b w(x)
w(x+h)$, we see that
\begin{align} \label{6.37}
\| \overline{w(x+h)b}\|^2=
\|\overline{w(x)b}\|^2+\|\overline{w(x)w(x+h)}\|^2 -2
\|\overline{w(x)b}\| \|\overline{w(x)w(x+h)}\| \cos\psi.
\end{align}


\begin{figure}
\psfig{figure=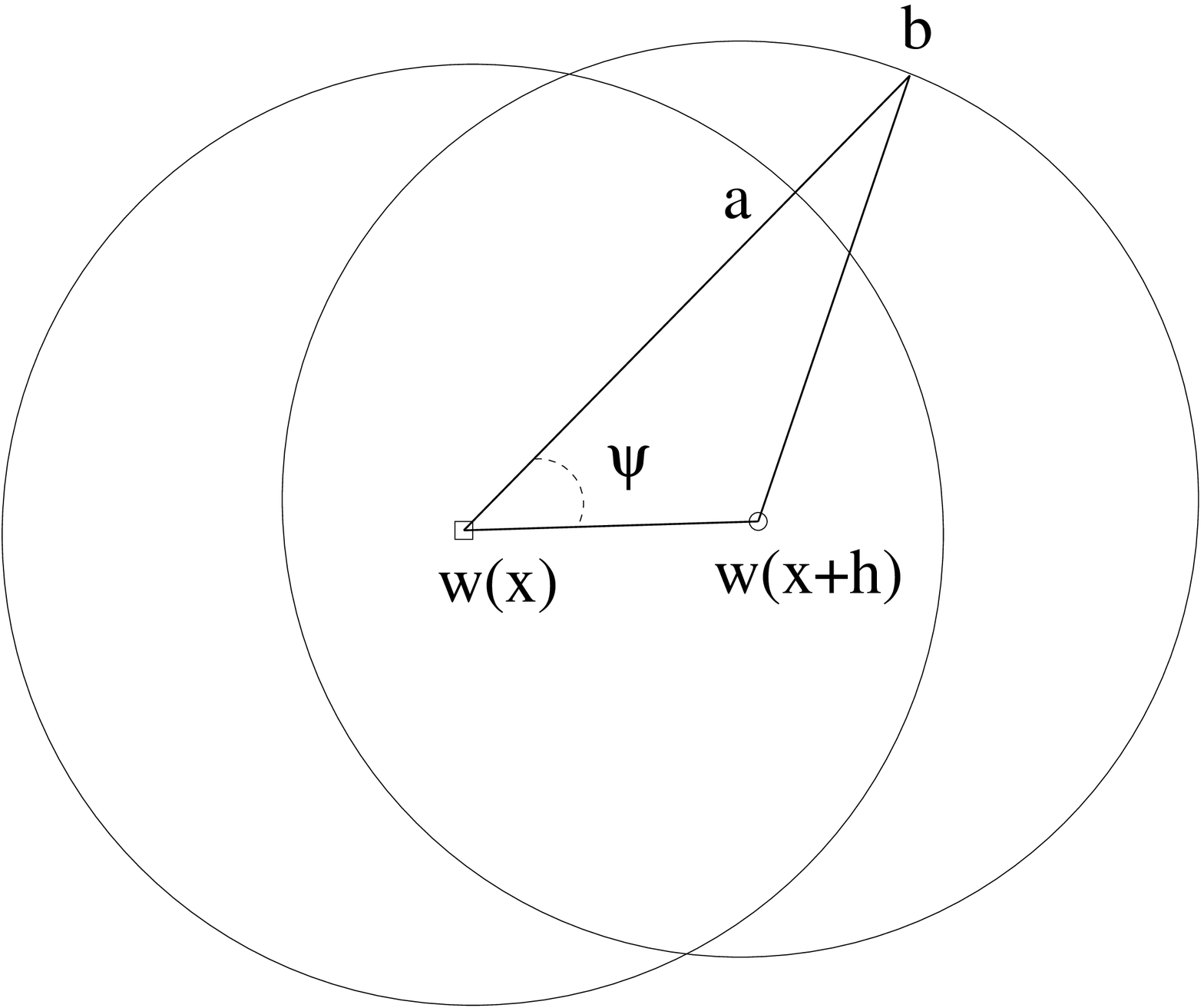, width=3in}
 \caption{  }
\end{figure}

Define $\|\overline{w(x)w(x+h)}\| = \rho(h)$ and
$\|\overline{w(x)b}\|=z$.
 Note that $\| \overline{w(x+h)b}\|= r(x+h)$. Then (\ref{6.37}) may be
 reformulated as a quadratic equation with respect to $z$:
$$
z^2-2\rho z \cos\psi+\rho^2 -r^2(x+h) = 0
$$
and therefore $z=\rho\cos\psi+\sqrt{r^2(x+h)-\rho^2\sin^2\psi}$.
We choose the positive sign because for $\psi=0$, we should get
$z= \rho+ r(x+h)$. As a result,
\begin{align} \label{6.38}
 \|\overline{ab}(\psi)\|=z-r(x)=\rho\cos\psi+\sqrt{r^2(x+h)-\rho^2\sin^2\psi}-r(x).
\end{align}
We introduce polar coordinates in the set (\ref{6.35}):
$$
\cD(x+h)\setminus \cD(x) \ni w \to (\psi, \om, \gamma)
$$
with $(\psi, \om)$ the spherical coordinate: $\psi$ is the angle
in Figure 1 and $\om$ is complement spherical coordinate; the
coordinate $\gamma$ is the distance from $a$ to the point that has
coordinate $w$ where $a$ is the following point:
$a=\overline{w(x)w} \cap \partial B(x)$. As is well-known, the
Jacobian of the transformation $w=w(\psi, \om, \gamma)$ is equal
to $(r(x)+\gamma)^{m-1}$. So
\begin{align} \label{6.39}
dw=(r(x)+\gamma)^{m-1} d\psi d\om d\gamma.
\end{align}

Let us estimate (\ref{6.38}) for small $\|h\|$. By (\ref{6.28}),
\begin{align} \label{6.40}
\rho^2(h)=\|(E+\a^{\ast}\a)^{-1}\a^{\ast}h\|^2+\|\a(E+\a^{\ast}\a)^{-1}\a^{\ast}h\|^2
\end{align}
and using (\ref{6.2}),
\begin{align*}
r^2(x+h)
&=r^2(x)-2((E+\a^{\ast}\a)^{-1}\a^{\ast}x,(E+\a^{\ast}\a)^{-1}\a^{\ast}h)
    \\
&-2(x-\a(E+\a^{\ast}\a)^{-1}\a^{\ast}x,
h-\a(E+\a^{\ast}\a)^{-1}\a^{\ast}h) + O^2(h).
\end{align*}

Hence by Taylor expansion and by (\ref{6.10'}),
\begin{align} \label{6.41}
\|\overline{ab}(\psi)\|(x,h)=\rho(h)\cos\psi-\frac1{r(x)}((E-\a(E+\a^{\ast}\a)^{-1}\a^{\ast})x,
h) + o(\|h\|^2).
\end{align}

Let $\pd_1 (x,h)$ be the part of the boundary for the set
$\cD(x+h)\setminus \cD(x)$, composed of the points on sphere $\pd
\cD(x)$. Changing to polar coordinates, and applying the Taylor
expansion for $\Gamma$ and using (\ref{6.39}) and (\ref{6.41}), we
get
\begin{align} \label{6.42}
 \int_{\cD(x+h)\setminus \cD(x)}\Gamma(w,x)dw
&=
\int_{\pd_1(x,h)}\int_0^{\|\overline{ab}(\psi)\|(x,h)}\Gamma(\psi,\om,\gamma,x)
(r(x)+\gamma)^{m-1}d\gamma d\psi d\om
\nonumber \\
& = \int_{\pd_1(x,h)}\Gamma(\psi,\om,0,x)\Psi_1(x,\psi,e) d\psi
d\om \|h\| + O(\|h\|^2),
\end{align}
where, recall, $e=h/\|h\|$, and the function $\Psi (x,\psi,e)$ is
defined by
\begin{align*}
\int_0^{\|\overline{ab}(\psi)\|(x,h)}(r(x)+\gamma)^{m-1}d\gamma
 & =
\frac1{m} ((r(x)+\|\overline{ab}(\psi)\|(x,h))^m-r(x)^m)  \nonumber \\
 & =
\Psi (x,\psi,e) \|h\| + o(\|h\|^2).
\end{align*}
In other words, by (\ref{6.41}),
\begin{align} \label{6.43}
\Psi
(x,\psi,e)=r^{m-1}(x)(\rho(e)\cos\psi-\frac1{r(x)}((E-\a(E+\a^{\ast}\a)^{-1}\a^{\ast})x,e)).
\end{align}


\begin{figure}
\psfig{figure=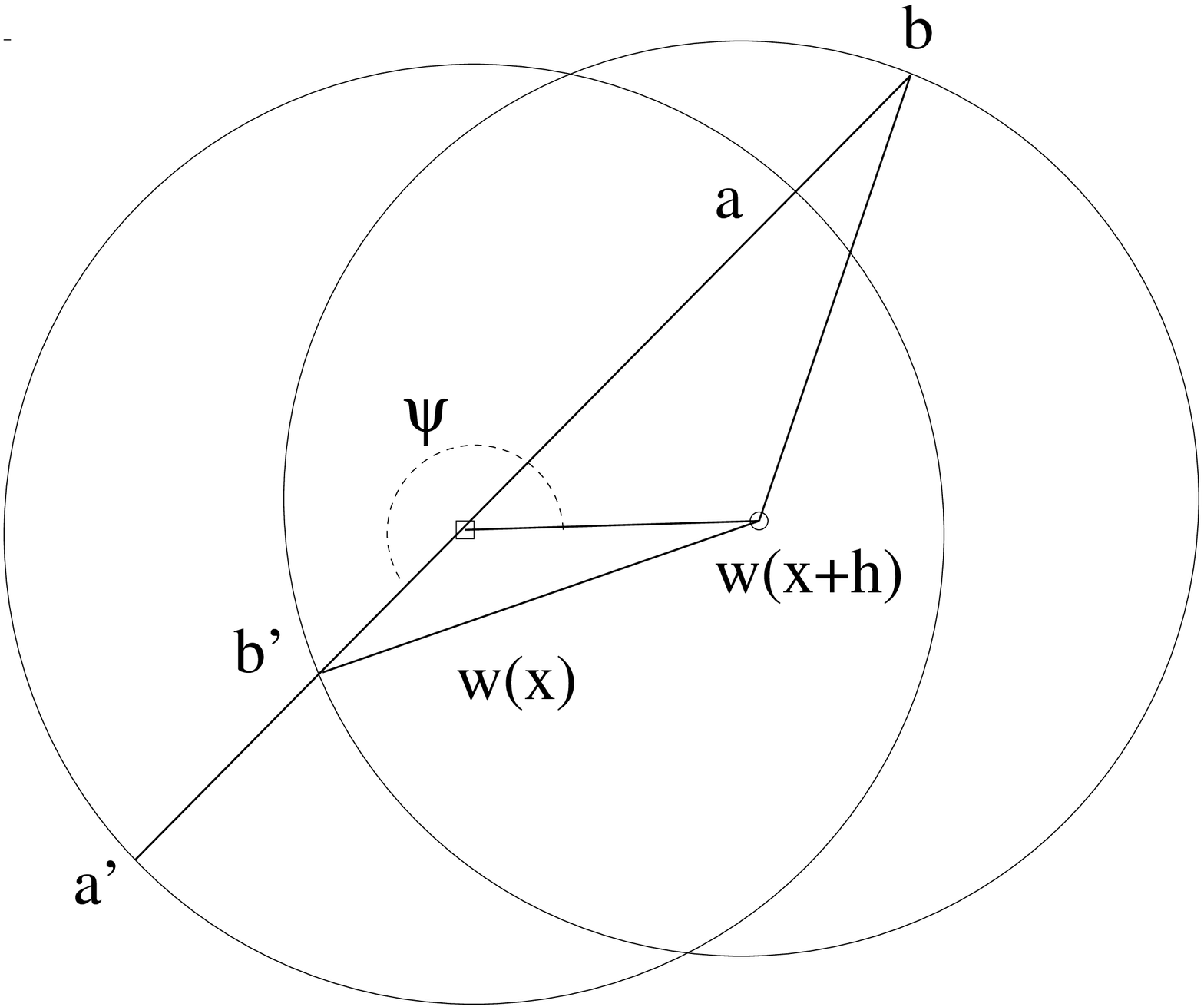,width=3.5in}
 \caption{  }
\end{figure}

To calculate the integral over $\cD(x)\setminus \cD(x+h)$ in
(\ref{6.33}) we use notations of points $w(x), w(x+h), a', b'$ and
angle $\psi $ from Figure 2. Similarly to the case
$\cD(x+h)\setminus \cD(x)$ we define $\|\overline{w(x)w(x+h)}\| =
\rho(h),\|\overline{w(x)b'}\|=z, \| \overline{w(x+h)b')}\|=
r(x+h)$. Then by Cosine Theorem in triangle $\bigtriangleup
b'w(x)w(x+h)$ we have the equality $z^2+\rho^2-2z\rho \cos{\psi}-
r^2(x+h)=0$ and therefore
$z=\rho\cos{\psi}+\sqrt{r^2(x+h)-\rho^2\sin^2{\psi}}$ (We put sign
``plus" before square root because for $\psi =\pi $ the equality
$z=r(x+h)-\rho $ holds). Then
$$
\| \overline{a'b'}\| =r(x)-z=-\rho (h)\cos{\psi}+
r^{-1}(x)((E-\alpha (E+\alpha^*\alpha )^{-1}\alpha^*)x,h)
$$
and similarly to (\ref{6.42}),(\ref{6.43}) we get
\begin{align} \label{6.44}
-\int_{\cD(x)\setminus \cD(x+h)}\Gamma(w,x)dw
 &=\int_{\pd_2(x,h)}\Gamma(\psi,\om,0,x)\Psi (x,\psi,e) d\psi d\om \|h\|+
o(\|h\|^2),
\end{align}
where $\Psi(x,\psi,e)$ is defined in (\ref{6.43}) and $\pd_2
(x,h)$ is the part of the boundary for the set $\cD(x)\setminus
\cD(x+h)$   composed of the points on sphere $\pd \cD(x)$.

Substituting (\ref{6.42}),(\ref{6.43}),(\ref{6.44}) into
(\ref{6.33}), dividing the obtained equation by $\l =\|h\|$ and
taking the limit $\l \downarrow 0$, we conclude
\begin{align*}
 \lim_{\lambda \downarrow 0} \frac{P(x+\lambda e)-P(x)}{\lambda}
 &=\int_{\pd \cD(x)}\Gamma(\psi,\om,0,x)\Psi(x,\psi,e) d\psi d\om +
 \int_{\cD(x)} (\Gamma_x'(w,x),e)dw,
\end{align*}
This equality shows that the density $P(x)$ possesses the first
variation \footnote{Recall (see, for instance, \cite{ATF}) that
$P(x)$ possesses the first variation at a point $x$ if for each
$h\in X_{\si \si_k}$ there exists a limit $\lim_{\lambda
\downarrow 0}(P(x+\lambda h)-P(x))/\lambda:=P'(x,h)$. Evidently,
$P'(x,h)$ is positively homogeneous on $h:\; P'(x,\lambda h)=
\lambda P'(x,h),\; \forall \lambda >0.$ The first variation
$P'(x,h)$ called Lagrange variation if $P'(x,-h)=-P(x,h)$.}
 for each $x$ on $\text{Int} (\pi B)$, and moreover
\begin{align} \label{6.45}
 P'(x,h)
  &=\int_{\pd \cD(x)}\Gamma(\psi,\om,0,x)\Psi(x,\psi,h) d\psi d\om +
 \int_{\cD(x)} (\Gamma_x'(w,x),h)dw.
\end{align}

\begin{theorem}\label{th7.1}
The first variation $P'(x,h)$ is the Lagrange variation. Moreover,
the function
\begin{equation}\label{7.10}
   x\rightarrow \sup_{\| e\|_{X_{\si \si_k}}} |P'(x,e)|\quad \mbox{is continuous for}\;
x\in \mbox{Int}(\pi B)
\end{equation}
\end{theorem}

\begin{proof} Note that the second terms in right sides of
(\ref{6.45}),(\ref{6.43}) are linear with respect to $h$ (to $e$).
Calculation of (\ref{6.43}) in the case of $-e$ gives that the
first term in right side of (\ref{6.43}) is equal $r^{m-1}(x)\rho
(-e)\cos{\psi_1}$ where $\psi_1=\psi +\pi$ and $\psi$ is the angle
from (\ref{6.43}). Hence,
$$
    r^{m-1}(x)\rho (-e)\cos{\psi_1}=r^{m-1}(x)\rho (e)\cos{(\psi +\pi )}=
-r^{m-1}(x)\rho (e)\cos{\psi}.
$$
Therefore $P'(x,-h)=-P'(x,h)$.

The assertion (\ref{7.10}) follows directly from the explicite
formulas (\ref{6.45}),(\ref{6.43}).
\end{proof}

Now we are in position to prove that $P(x)$ satisfy (\ref{3.201}).

\begin{lem}\label{lem7.1}
For every $v_1,v_2\in X_{\si \si_{k}}$
\begin{eqnarray}\label{7.1}
  \int_{X_{\si \si_k}}|P(x-v_1)-P(x-v_2)|\, dx\le c\| v_1-v_2\|_{X_{\si \si_k}}
\end{eqnarray}
where $c>0$ does not depend on $v_1,v_2$.
\end{lem}

\begin{proof}
For each $x,v_1,v_2\in X_{\si \si_k}$
\begin{equation}\label{7.2}
    P(x-v_1)-P(x-v_2)=\int_0^1\frac{dP(x+v_2+\theta (v_2-v_1))}{d\theta }\,d\theta
\end{equation}
Note that derivative $dP/d\theta$ is well defined for every
$\theta \in [0,1]$ except, maybe, one value such that
$x+v_2+\theta (v_2-v_1)\in \partial \pi B$ because $P=0$ outside
$\pi B$ and $P$ possesses Lagrange variation inside $\pi B$.
Besides, by Theorems \ref{th7.1} and \ref{theorem6.1} the function
$\theta \rightarrow dP/d\theta $ is integrable and, hence,
(\ref{7.2}) is well defined. Therefore (\ref{7.2}) implies
$$
 \int_{X_{\si si_k}}|P(x-v_1)-P(x-v_2)|\,dx
$$
$$
\le \| v_2-v_1\|_{X_{\si \si_k}}\int_{X_{\si \si_k}}\int_0^1
|P'(x+v_2+\theta (v_2-v_1),\frac{v_2-v_1}{\| v_2-v_1\|_{X_{\si
\si_k}}})|\,d\theta dx
$$
$$
\le \| v_2-v_1\|_{X_{\si \si_k}}\int_{X_{\si
\si_k}}\sup_{\|e\|_{X_{\si \si_k}}=1}|P'(x,e)|\,dx
$$
This implies inequality (\ref{7.1})
\end{proof}

Thus, we have proved that RDS (\ref{3.1}) is a particular case of
RDS (\ref{3.18}) and therefore it satisfies all conditions of
Theorem \ref{kuksin}. Indeed, Lemma \ref{lem7.1} implies that
random dynamical system (\ref{3.1}) satisfies the condition
(\ref{3.201}). As was shown above, conditions
(\ref{3.19}),(\ref{3.20}) are true for RDS (\ref{3.1}) in virtue
of Lemma \ref{lem4.10}. So assertion of Theorem \ref{kuksin} is
true for RDS (\ref{3.1}). This proves Theorem \ref{main} and
completes our investigation in this paper.


\bigskip
\bigskip
\bigskip

{\bf Acknowledgement.}  A part of this work was done   while A. V.
Fursikov was visiting Department of Applied Mathematics, Illinois
Institute of Technology, Chicago, USA. This visit was sponsored by
the National Research Council.


\end{document}